%
\def\conv{\mathop{\vrule height2,6pt depth-2,3pt 
    width 5pt\kern-1pt\rightharpoonup}}

\advance\vsize by 1 true cm
%
\def\dess #1 by #2 (#3){
  \vbox to #2{
    \hrule width #1 height 0pt depth 0pt
    \vfill
    \special{picture #3} 
    }
  }

\def\dessin #1 by #2 (#3 scaled #4){{
  \dimen0=#1 \dimen1=#2
  \divide\dimen0 by 1000 \multiply\dimen0 by #4
  \divide\dimen1 by 1000 \multiply\dimen1 by #4
  \dess \dimen0 by \dimen1 (#3 scaled #4)}
  }
%
\def \trait (#1) (#2) (#3){\vrule width #1pt height #2pt depth #3pt}
\def \fin{\hfill
	\trait (0.1) (5) (0)
	\trait (5) (0.1) (0)
	\kern-5pt
	\trait (5) (5) (-4.9)
	\trait (0.1) (5) (0)
\medskip}
%


\font\sevenbf=cmbx7

\baselineskip=15pt
\abovedisplayskip=15pt plus 4pt minus 9pt
\belowdisplayskip=15pt plus 4pt minus 9pt
\abovedisplayshortskip=3pt plus 4pt
\belowdisplayshortskip=9pt plus 4pt minus 4pt
\let\epsilon=\varepsilon

\def\biblio #1 #2\par{\parindent=30pt\item{}\kern -30pt\rlap{[#1]}\kern
30pt #2\smallskip}
 %
\catcode`\@=11
\def\@lign{\tabskip=0pt\everycr={}}
\def\equations#1{\vcenter{\openup1\jot\displ@y\halign{\hfill\hbox
{$\@lign\displaystyle##$}\hfill\crcr
#1\crcr}}}
\catcode`\@=12
%
\def\pmb#1{\setbox0=\hbox{#1}%
\hbox{\kern-.04em\copy0\kern-\wd0
\kern.08em\copy0\kern-\wd0
\kern-.02em\copy0\kern-\wd0
\kern-.02em\copy0\kern-\wd0
\kern-.02em\box0\kern-\wd0
\kern.02em}}
%
\def\undertilde#1{\setbox0=\hbox{$#1$}
\setbox1=\hbox to \wd0{$\hss\mathchar"0365\hss$}\ht1=0pt\dp1=0pt
\lower\dp0\vbox{\copy0\nointerlineskip\hbox{\lower8pt\copy1}}}
%

%

\def\maj#1#2,{\rm #1\sevenrm #2\rm{}}
\def\Maj#1#2,{\bf #1\sevenbf #2\rm{}}
\outer\def\lemme#1#2 #3. #4\par{\medbreak
\noindent\maj{#1}{#2},\ #3.\enspace{\sl#4}\par
\ifdim\lastskip<\medskipamount\removelastskip\penalty55\medskip\fi}

\def\Remark #1. {\noindent{\Maj REMARK,\ \bf #1. }}

\outer\def\Lemme#1#2 #3. #4\par{\medbreak
\noindent\Maj{#1}{#2},\ \bf #3.\rm\enspace{\sl#4}\par
\ifdim\lastskip<\medskipamount\removelastskip\penalty55\medskip\fi}



\def\Notation #1. {\noindent{\Maj NOTATION,\ \bf #1. }}

\def\Example #1. {\noindent{\Maj EXAMPLE,\ \bf #1. }}

\hfuzz=1cm


\catcode`\ˆ=\active     \def ˆ{\`a}
\catcode`\‰=\active     \def ‰{\^a}
\catcode`\=\active     \def {\c c}
\catcode`\Ž=\active    \def Ž{\'e} 

\catcode`\=\active   \def {\`e}
\catcode`\=\active   \def {\^e}
\catcode`\'=\active   \def '{\"e}
\catcode`\"=\active   \def "{\^\i}
\catcode`\•=\active   \def •{\"\i}
\catcode`\™=\active   \def ™{\^o}
\catcode`\š=\active   \defš{}
\catcode`\=\active   \def {\`u}
\catcode`\ž=\active   \def ž{\^u}
\catcode`\Ÿ=\active   \def Ÿ{\"u}
\catcode`\ =\active   \def  {\tau}
\catcode`\¡=\active   \def ¡{\circ}
\catcode`\¢=\active   \def ¢{\Gamma}
\catcode`\¤=\active   \def ¤{\S\kern 2pt}
\catcode`\¥=\active   \def ¥{\puce}
\catcode`\§=\active   \def §{\beta}
\catcode`\¨=\active   \def ¨{\rho}
\catcode`\©=\active   \def ©{\gamma}
\catcode`\­=\active   \def ­{\neq}
\catcode`\°=\active   \def °{\ifmmode\ldots\else\dots\fi}
\catcode`\±=\active   \def ±{\pm}
\catcode`\²=\active   \def ²{\le}
\catcode`\³=\active   \def ³{\ge}
\catcode`\µ=\active   \def µ{\mu}
\catcode`\¶=\active   \def ¶{\delta}
\catcode`\·=\active   \def ·{\Sigma}
\catcode`\¸=\active   \def ¸{\Pi}
\catcode`\¹=\active   \def ¹{\pi}
\catcode`\»=\active   \def »{\Upsilon}
\catcode`\¾=\active   \def ¾{\alpha}
\catcode`\À=\active   \def À{\cdots}
\catcode`\Â=\active   \def Â{\lambda}
\catcode`\Ã=\active   \def Ã{\sqrt}
\catcode`\Ä=\active   \def Ä{\varphi}
\catcode`\Å=\active   \def Å{\xi}
\catcode`\Æ=\active   \def Æ{\Delta}
\catcode`\Ç=\active   \def Ç{\cup}
\catcode`\È=\active   \def È{\cap}
\catcode`\Ï=\active   \def Ï{\oe}
\catcode`\Ñ=\active   \def Ñ{\to}
\catcode`\Ò=\active   \def Ò{\in}
\catcode`\Ô=\active   \def Ô{\subset}
\catcode`\Õ=\active   \def Õ{\superset}
\catcode`\Ö=\active   \def Ö{\over}
\catcode`\×=\active   \def ×{\nu}
\catcode`\Ù=\active   \def Ù{\Psi}
\catcode`\Ú=\active   \def Ú{\Xi}
\catcode`\Ü=\active   \def Ü{\omega}
\catcode`\Ý=\active   \def Ý{\Omega}
\catcode`\ß=\active   \def ß{\equiv}
\catcode`\à=\active   \def à{\chi}
\catcode`\á=\active   \def á{\Phi}
\catcode`\ä=\active   \def ä{\infty}
\catcode`\å=\active   \def å{\zeta}
\catcode`\æ=\active   \def æ{\varepsilon}
\catcode`\è=\active   \def è{\Lambda}  
\catcode`\é=\active   \def é{\kappa}
\catcode`\ë=\active   \defë{\Theta}
\catcode`\ì=\active   \defì{\eta}
\catcode`\í=\active   \defí{\theta}
\catcode`\î=\active   \defî{\times}
\catcode`\ñ=\active   \defñ{\sigma}
\catcode`\ò=\active   \defò{\psi}
\def\date{\number\day\
\ifcase\month \or janvier \or f\'evrier \or mars \or avril \or mai \or juin \or juillet \or ao\^ut  \or
septembre \or octobre \or novembre \or d\'ecembre \fi
\ \number\year}

\font \Ggras=cmb10 at 12pt
\font \ggras=cmb10 at 11pt

\def\GS{{\bf \Ggras S}}

\def\liminf{\mathop{\underline{\rm lim}}}
\def\limsup{\mathop{\overline{\hbox{\rm lim}}}}
\def\sym{\fam\comfam\com}
\font\tensym=msbm10
\font\sevensym=msbm7
\font\fivesym=msbm5
\newfam\symfam
\textfont\symfam=\tensym
\scriptfont\symfam=\sevensym
\scriptfont\symfam=\fivesym
\def\sym{\fam\symfam\relax}
\def\N{{\sym N}}
\def\R{{\sym R}}

\def\Z{{\sym Z}}


\centerline{\Ggras  Asymptotic behavior of Structures made of Plates}
\vskip 5mm 
\centerline{ G. Griso}
\centerline{  Laboratoire J.-L. Lions--CNRS, Bo\^\i te courrier 187, }
\centerline{Universit\'e  Pierre et Marie Curie, 4~place Jussieu, 75005 Paris, France,}
\centerline{ Email: griso@ann.jussieu.fr}
\vskip 8mm

\noindent{\ggras \bf Abstract. } {\sevenrm   The aim of this work is to study the asymptotic behavior of a structure made of
plates of thickness $\scriptstyle 2\delta$   when $\scriptstyle\delta\to 0$. This study is carried on within the frame of linear
elasticity by using the unfolding method. It is based on several decompositions of the structure displacements and on the passing
to the limit in fixed domains.

We begin with studying  the displacements of a plate. We show that any displacement is the sum of an elementary displacement
concerning the normal lines on the middle surface of the plate and a residual displacement linked to these normal lines
deformations. An elementary displacement is linear with respect to the variable $\scriptstyle x_3$. It is written $\scriptstyle {\cal
U}(\hat x)+{\cal R}(\hat x)\land x_3\vec e_3$ where  $\scriptstyle{\cal U}$ is a displacement of the mid-surface of the plate.
We show a priori estimates and convergence results when $\scriptstyle\delta\to 0$. We characterize the limits of the unfolded
displacements of a plate as well as the limits of the unfolded of the strained tensor.

Then we extend these results to the structures made of plates. We show that any displacement of a structure is the sum of an
elementary displacement  of each plate and of a residual displacement. The elementary displacements of the structure  (e.d.p.s.)
 coincide with  elementary rods displacements in the junctions.  Any  e.d.p.s. is given by two functions belonging to $\scriptstyle
H^1(\GS;\R^3)$ where $\scriptstyle \GS$ is the skeleton of the structure (the plates mid-surfaces set). One of these functions :
$\scriptstyle {\cal U}$ is the skeleton displacement. We show that $\scriptstyle {\cal U}$ is the sum of an extensional
displacement and of an inextensional one. The first one characterizes the membrane displacements and the second one is a rigid
displacement in the direction of the plates and it characterizes the plates flexion.

Eventually we pass to the limit as $\scriptstyle \delta\to 0$ in the linearized elasticity system,   on the one
hand we obtain   a variational  problem that is satisfied by the  limit extensional displacement,  and on the other hand,  a  
variational  problem  satisfied by  the limit of   inextensional  displacements.  } 
\smallskip
\noindent{\ggras \bf R\'esum\'e. } {\sevenrm  Le but de ce travail est d'\'etudier le comportement asymptotique d'une
structure formŽe de plaques d'Žpaisseur $\scriptstyle 2\delta$  lorsque $\scriptstyle\delta\to 0$. Cette \'etude est men\'ee dans
le cadre de l'\'elasticit\'e lin\'eaire en utilisant la m\'ethode de l'\'eclatement.   Elle est bas\'ee sur plusieurs d\'ecompositions des
d\'eplacements de la structure, et  sur le passage \`a la limite dans des domaines fixes.  

On commence par une Žtude des dŽplacements d'une plaque. On montre que tout d\'eplacement d'une 
plaque est la somme d'un d\'eplacement  \'el\'ementaire  concernant  les normales ˆ la surface moyenne de la plaque  et d'un 
d\'eplacement  r\'esiduel  li\'e aux  d\'eformations  de ces normales. Un d\'eplacement  \'el\'ementaire est affine par rapport
\`a la variable $\scriptstyle x_3$, il s'Žcrit $\scriptstyle {\cal U}(\hat x)+{\cal R}(\hat x)\land x_3\vec e_3$ o
$\scriptstyle{\cal U}$ est un dŽplacement de la surface moyenne de la plaque. On \'etablit des estimations a  priori   et  des
r\'esultats de convergence  lorsque  $\scriptstyle\delta\to 0$.  On caract\'erise  les limites des \'eclat\'es des  d\'eplacements
d'une plaque,  ainsi que les limites des \'eclat\'es du tenseur des d\'eformations.   

On Žtend ensuite ces rŽsultats aux structures formŽes de plaques. On montre que tout dŽplacement d'une structure est la somme
d'un dŽplacement ŽlŽmentaire de chaque plaque et d'un dŽplacement rŽsiduel. Les dŽplacements ŽlŽmentaires de la structure
(d.e.s.p.)  co\"\i  ncident avec des dŽplacements \'el\'ementaires de poutres dans les jonctions.  Tout d.e.s.p.  est donnŽ par deux
fonctions  appartenant ˆ $\scriptstyle H^1(\GS;\R^3)$ o $\scriptstyle \GS$ est le squelette de la structure (l'ensemble des
surfaces moyennes  des plaques). L'une de ces fonctions : $\scriptstyle {\cal U}$ est le dŽplacement du squelette.  
On montre  que $\scriptstyle {\cal U}$   est la  somme  d'un d\'eplacement  extensionnel  et d'un d\'eplacement   inextensionnel.
Le  premier caractŽrise les dŽplacements membranaires des surfaces moyennes, le second est un dŽplacement rigide dans la
direction  des plaques; il caractŽrise la flexion des plaques. 

 Pour finir on passe  \`a la  limite  pour $\scriptstyle\delta\to 0$  dans le syst\`eme de l'\'elasticit\'e lin\'eaire, on obtient d'une
part  un   probl\`eme variationnel  v\'erifi\'e par la limite  des  d\'eplacements  extensionnels,   et d'autre part  un   probl\`eme
variationnel  v\'erifi\'e par  la limite  des  d\'eplacements    inextensionnels.}

\vfill\eject
\noindent {\ggras  1. Introduction} 
\vskip 6pt 
Many articles and books have been dedicated to the mathematical justification of plates models (see for example [1,2]).  A first
study concerning the asymptotic behavior of a structure made of two thin plates of thickness
$\varepsilon$,  is due to Le Dret [7].  The  obtained asymptotic model derives from the  three-dimensional  system of
elasticity thanks to a  thin domain standard technique (the  plates are transformed into a fixed domain). At the limit, Le Dret
obtains a two-dimensional system coupling the  flexion displacements of the two mid-surfaces of the plates.  

Our study  continues [4] and [5]. In this paper we use again the  notions of elementary displacements and of  extensional and
inextensional displacements and we extend them to the plates displacements and to the displacements of structures made of
plates.  Our paper is organised into three parts. In the first one we study the displacements of a plate, the second one  is devoted
to the displacements of a structure made of plates from which we deduce  the asymptotic behavior of a structure made of thin
plates. And in the third part  we prove the technical lemmas used in the two first parts of our paper.

In Section 2 we consider a plate of  thickness $2\delta$. We first introduce the elementary displacements of a plate (Definition
2.1).  These  are the displacements of the normal lines of the mid-surface of the plate.  An elementary displacement is linear with
respect to the variable $x_3$. It is written ${\cal U}(\hat x)+{\cal R}(\hat x)\land x_3\vec e_3$ where  ${\cal U}$ is a
displacement of the mid-surface. By such a displacement the normal line  is transformed into a line which is  generally no longer
perpendicular to the mid-surface.  With each displacement $u$ of the plate  we associate an elementary displacement $U_e$
(Definition 2.2).  Theorem 2.3  gives estimates of appropriate norms of $U_e$ and of the displacement $u- U_e$ in terms of 
$\delta$. Using the elementary displacement $U_e$ we show (formula (2.3)) that the displacement $u$ is the sum of a
Kirchhoff-Love displacement and of a residual  one  $\widetilde{u}$, which satisfies estimate (2.4).  We are now equipped to obtain
the asymptotic behavior of a displacements  sequence  $\bigl( u_\delta\bigr)_{\delta>0}$ with strain energy of order $\delta$.
This is the main result of this section  and it is given in Theorem 2.6. The previous decomposition allows us to give a simple
interpretation (see Theorem 2.6)  of the limits of  the unfolding ${\cal T}_\delta\bigl(\gamma_{ij}(u_\delta)\bigr)$ of the strain
tensor $\gamma_{ij}(u_\delta)$ (where the unfolding operator ${\cal T}_\delta$ is given in Definition 2.5) in terms of the
derivatives limits of the Kirchhoff-Love displacements and of  the residual  displacements.  There is not a unique
associated elementary displacement that satisfies estimates $(2.2)$. In Definition 2.2 we give the simplest one. But the one we
give in Definition 2.9 is more suitable for the study of  a structure made of plates.

The structure ${\cal S}_\delta$  made of plates of thickness $2\delta$ is  introduced in Section  3. Our hypotheses about the
skeleton of the structure $\GS$  (i.e. the plates mid-surfaces set) allow us to consider a wide range of structures. We
extend to them the notions and decompositions of Section 2.  Definition 3.1 gives us the elementary displacements of
plates-structure  (e.d.p.s.). These displacements coincide with elementary plate displacements in each plate and there are rods
elementary displacements in the junctions (see [5]).  Any  e.d.p.s. is known by two functions belonging to $H^1(\GS;\R^3)$. The
first one ${\cal U}$ is the skeleton displacement, the second one gives the rotations of the normal lines of the mid-surfaces. We
show that
${\cal U}$ is the sum of an extensional displacement and of an inextensional one  (Definitions 3.6 and 3.5). The first one
characterizes the membrane displacements and the second one is a rigid displacement in the direction of the plates and it
characterizes the plates flexion. Corollary of Lemma 3.7 gives estimates for them with an appropriate norm. In subsection 3.4. we
consider an e.d.p.s. sequence $\bigl( u_\delta\bigr)_{\delta>0}$  with strain energy of order $\delta$. Thanks to all these
decompositions  we give the limits of  the unfolding
${\cal T}_\delta\bigl(\gamma_{ij}(u_\delta)\bigr)$ of the strain tensor as in the case of a plate. We also characterize the space
of the inextensional limits displacements. In the last subsection, we  give  the limit for $\delta\to 0$ of  the  linearized  elasticity
system   (3.9), written in
${\cal S}_\delta$, where the applied forces $F_\delta$ satisfy assumptions (3.11).  The main results are Theorem 3.8 and  Theorem
3.10. In the first one we show  that the extensional displacement limit is the solution of  a second-order system, and in the second
one  we show that  the  limit  of the inextensional displacement is the solution of a fourth-order system. 
\vskip 1mm 
In  this work we use the Einstein convention of summation over repeated indices. As a rule, the Greek indices
$\alpha$ and $\beta$ take values in $\{1, 2\}$ and the Latin indices $i$, $i^{'}$, $j$ and $j^{'}$  take values in $\{1, 2, 3\}$. 
\bigskip
\noindent {\bf  2. The plate displacements} 
\vskip 1mm
\noindent {\bf 2.1. The elementary plate displacements} 
\vskip 1mm
\noindent  The Euclidian space $\R^3$ is related to the frame  $(O;\vec e_1,\vec e_2,\vec e_3)$.  Let
$\omega$ be a bounded domain in $\R^2$ with a lipschitzian boundary. The plate
$\Omega_\delta=\omega\times]-\delta,\delta[$,
$\delta>0$, is the open set having as middle surface $\omega$ and as thickness $2\delta$. The direction of the normal lines of
$\omega$ is given by $\vec e_3$.The reference plate is the open set $\Omega=\omega\times]-1,1[$.

\noindent The running point of $\Omega_\delta$ (respectively $\Omega$) is denoted $x=(x_1,x_2,x_3)=(\widehat{ x},x_3)$, 
(resp.  $(\widehat{x}, t_3)$) where $\widehat{x}\in \omega$ and $t_3\in ]-1,1[$. 

\noindent For any open set $\omega^{'}$ of $\R^n$, $n\in\{2,3\}$, and any displacement $u$ belonging to
$H^1(\omega^{'},\R^n)$, we put
$${\cal E}(u,\omega^{'})=\int_{\omega^{'}}\gamma_{ij}(u)\gamma_{ij}(u),\qquad \gamma_{ij}(u)={1\over 2}\Bigl\{{\partial
u_i\over \partial x_j}+{\partial u_j\over \partial x_i}\Bigr\}, \qquad {\cal D}(u,\omega^{'})=\int_{\omega^{'}}{\partial u_i\over
\partial x_j}{\partial u_i\over \partial x_j}.$$
\noindent{\bf Definition  2.1 : } An {\bf elementary plate displacement} (e.p.d.) is an element $\Phi$ belonging to
$H^1(\Omega_\delta,\R^3)$, such that 
$$\Phi(x)={\cal A}(\widehat{ x})+{\cal B}(\widehat{ x})\land x_3\vec e_3,\qquad\hbox{a.e.}\enskip x=(\widehat{ x},x_3)
\in\omega\times ]-\delta,\delta[=\Omega_\delta,$$ where ${\cal A}$ and ${\cal B}$ belong to $H^1(\omega ,\R^3)$;
${\cal A}$ is the first component and ${\cal B}$  the second component of the e.p.d. $\Phi$.
\vskip 1mm
\noindent{\bf Elementary plate displacement associated with a displacement of
$H^1(\Omega_\delta,\R^3)$.}
\vskip 1mm
\noindent{\bf Definition 2.2 : } With any displacement $u\in H^1(\Omega_\delta,\R^3)$, we associate the elementary plate
displacement $U_e$ defined as
$$\left\{\eqalign{  
U_e(x)&={\cal U}(\widehat{ x})+{\cal R}(\widehat{ x})\land x_3\vec e_3,\qquad x=(\widehat{ x},x_3)\in
\Omega_\delta,\cr {\cal U}(\widehat{ x})&={1\over 2\delta}\int_{-\delta}^\delta u(\widehat{ x},x_3)dx_3,\qquad
{\cal R}(\widehat{ x})={3\over 2\delta^3}\int_{-\delta}^\delta x_3\vec e_3\land u(\widehat{x},x_3) dx_3.\cr  }\right.\leqno
(2.1)$$ The component ${\cal R}_3$ of ${\cal R}$ is equal to 0.

\noindent {\bf Theorem 2.3  : }{\it The elementary plate displacement $U_e$ verifies 
$${\cal E}(U_e,\Omega_\delta)+{\cal D}(u-U_e,\Omega_\delta)+{1\over \delta^2}||u-U_e||^2_{L^2(
\Omega_\delta, \R^3)}\le C{\cal E}(u,\Omega_\delta).\leqno(2.2)$$ The constants depend  only on $\omega$.}  
\vskip 1mm
\noindent{\bf Proof : } See Annex A.\fin
\noindent{\bf Proposition  2.4  : }{\it Any displacement $u$ belonging to $H^1(\Omega_\delta,\R^3)$ is the sum
of a Kirchhoff-Love displacement and a residual one $\widetilde {u}$
$$u(x)=\Bigl({\cal U}_1(\widehat{x})-x_3{\partial {\cal U}_3\over \partial x_1}(\widehat{x})\Bigr)\vec e_1+\Bigl({\cal
U}_2(\widehat{x})-x_3{\partial {\cal U}_3\over \partial x_2}(\widehat{x})\Bigr)\vec e_2+ {\cal U}_3(\widehat{x})\vec
e_3+\widetilde{u}(x),\qquad x\in\Omega_\delta,\leqno (2.3)$$ where ${\cal U}$ is the first component of the e.p.d.
$U_e$. The residual displacement $\widetilde{u}$ belongs to $L^2(\omega,H^1(]-\delta,\delta[,\R^3))$ and verifies  
$${1\over \delta^2}||\widetilde{u}||^2_{L^2(\Omega_\delta,\R^3)}+\Bigl\|{\partial \widetilde{u}\over\partial
x_3}\Bigr\|^2_{L^2(\Omega_\delta,\R^3)}\le C{\cal E}(u,\Omega_\delta).\leqno (2.4)$$}
\noindent{\bf Proof : } We define $\widetilde{u}$ by
$$ \widetilde{u}(x)=u(x)-\Bigl({\cal U}_1(\widehat{x}) -x_3{\partial {\cal U}_3\over \partial x_1}(\widehat{x})\Bigr)\vec
e_1-\Bigl({\cal U}_2(\widehat{x})- x_3{\partial {\cal U}_3\over \partial x_2}(\widehat{x})\Bigr)\vec e_2-{\cal
U}_3(\widehat{x})\vec e_3.$$ This displacement belongs to 
$L^2(\omega, H^1(]-\delta,\delta[,\R^3))$.  We obtain $ (2.4)$ using the estimates of Theorem 2.3.\fin 
\noindent{\bf 2.2. Limit of  a plate displacements sequence.}
\vskip 1mm
\noindent{\bf Definition 2.5 : }  The unfolding operator ${\cal T}_\delta$ from
$L^2(\Omega_\delta,\R^n)$ into
$L^2(\Omega,\R^n)$ is defined by
$${\cal T}_\delta(\phi)(\widehat{x},t_3)=\phi(\widehat{x},\delta t_3),\qquad \hbox{a.e. in }\; \Omega.\leqno(2.5)$$ For any
element $\phi\in H^1(\Omega_\delta)$, we have ${\cal T}_\delta(\phi)\in H^1(\Omega)$ and
$${\partial {\cal T}_\delta(\phi)\over \partial x_\alpha}={\cal T}_\delta\Bigl({\partial \phi\over \partial x_\alpha}\Bigr),\qquad 
{\partial {\cal T}_\delta(\phi)\over \partial t_3}=\delta{\cal T}_\delta\Bigl({\partial \phi\over \partial x_3}\Bigr).\leqno (2.6)$$

\noindent{\bf Theorem 2.6 : }{\it Let $\bigl(u_\delta\bigr)_{\delta>0}$ be a sequence of displacements of
$H^1(\Omega_\delta,\R^3)$ verifying 
$${\cal E}(u_\delta,\Omega_\delta)\le C\delta.\leqno (2.7)$$ There exist $(a_\delta, b_\delta)\in \R^3\times\R^3$
 and extracted sequences (still denoted in the same way), such that
$$\left\{\eqalign{
&{\cal U}_{1,\delta}-a_{1,\delta}+x_2 b_{3,\delta}\rightharpoonup U_1,\qquad {\cal
U}_{2,\delta}-a_{2,\delta}-x_1 b_{3,\delta}\rightharpoonup U_2\quad \quad\hbox{weakly in}\quad H^1(\omega)\cr
&  \delta\bigl\{{\cal U}_{3,\delta}-a_{3,\delta}+x_1b_{2,\delta}-x_2b_{1,\delta}\bigr\}\rightharpoonup U_3 
\quad\hbox{weakly in}\quad H^1(\omega).\cr }\right.\leqno (2.8)$$
\noindent  Moreover $U_3$ belongs to $H^2(\omega)$. We have the following weak convergences of the unfolded of
$u_\delta$,  $\widetilde{u}_\delta$ and of the components of the strain tensor :
$$\left\{\eqalign{ 
&{\cal T}_\delta\bigl(u_{1,\delta}-a_{1,\delta}+x_2 b_{3,\delta}\bigr)\rightharpoonup U_1-t_3{\partial
U_3\over \partial x_1},\qquad{\cal T}_\delta\bigl(u_{2,\delta}-a_{2,\delta}-x_1 b_{3,\delta}\bigr)\rightharpoonup
U_2-t_3{\partial U_3\over \partial x_2},\cr 
& \delta {\cal T}_\delta\bigl(u_{3,\delta}-a_{3,\delta}+x_1b_{2,\delta}-x_2b_{1,\delta}\bigr)\rightharpoonup
U_3,\quad \hbox{weakly in}\quad H^1(\Omega),\cr 
&{1\over \delta}{\cal T}_\delta\bigl(\widetilde{u}_\delta\bigr)\rightharpoonup
\widetilde{u}\qquad \hbox{weakly in}\quad L^2(\omega, H^1(]-1,1[,\R^3)),\cr &{\cal
T}_\delta\bigl(\gamma_{\alpha\beta}(u_\delta)\bigr)\rightharpoonup {1\over 2}\Bigl\{{\partial U_\alpha\over \partial
x_\beta}+{\partial U_\beta\over \partial x_\alpha}\Bigr\}-t_3{\partial^2U_3\over \partial x_\alpha\partial x_\beta},\quad 
{\cal T}_\delta\bigl(\gamma_{\alpha 3}(u_\delta)\bigr)\rightharpoonup {1\over 2}{\partial\widetilde{u}_\alpha\over \partial
t_3},\cr 
&{\cal T}_\delta\bigl(\gamma_{33}(u_\delta)\bigr)\rightharpoonup {\partial \widetilde{u}_3\over
\partial t_3}\qquad \hbox{weakly in}\quad L^2(\Omega).\cr}\right.\leqno (2.9)$$}
\vskip-0.4cm
\noindent{\bf Proof : } With each $u_\delta$ we associate the e.p.d.  $U_{e,\delta}$ with components
${\cal U}_\delta$ and ${\cal R}_\delta$. From $(2.4)$ the displacement ${\cal U}_{M,\delta}={\cal U}_{1,\delta}\vec e_1+{\cal
U}_{2,\delta}\vec e_2$   has a strain energy $\displaystyle {\cal E}({\cal U}_{M,\delta},\omega)\le
{C\over\delta}{\cal E}(U_{e,\delta},\Omega_\delta)\le {C\over\delta}{\cal E}(u_\delta,\Omega_\delta)\le C$.
The classical Korn inequality applied to ${\cal U}_{M,\delta}$ affirms the existence of a rigid displacement $\displaystyle
r_{M,\delta}(\widehat{x})=\pmatrix{a_{1,\delta}-b_{3,\delta}x_2\cr a_{2,\delta}+b_{3,\delta}x_1\cr}$, such that
$$||{\cal U}_{M,\delta}-r_{M,\delta}||^2_{L^2(\omega,\R^2)}+{\cal D}({\cal U}_{M,\delta}-r_{M,\delta},\omega)\le
C\displaystyle {\cal E}({\cal U}_{M,\delta},\omega)\le {C\over \delta}{\cal E}(u_\delta,\Omega_\delta)\le C.$$ If 
$b_{\alpha,\delta}$ is the mean of ${\cal R}_{\alpha,\delta}$ on $\omega$, we obtain from the
Poincar\'e-Wirtinger inequality
$$||{\cal R}_{\alpha,\delta}-b_{\alpha,\delta}||^2_{L^2(\omega)}\le \displaystyle C||\nabla {\cal
R}_{\alpha,\delta}||^2_{[L^2(\omega)]^2}\le {C\over\delta^3}{\cal E}(u_\delta,\Omega_\delta)\le {C\over\delta^2}$$ 
The estimate of ${\cal E}(U_{e,\delta},\Omega_\delta)$  obtained in Theorem 2.3, gives 
$$\Bigl\|{\partial {\cal U}_{3,\delta}\over \partial x_1}+{\cal R}_{2,\delta}\Bigr\|^2_{L^2(\omega)}+\Bigl\|{\partial {\cal
U}_{3,\delta}\over \partial x_2}-{\cal R}_{1,\delta}\Bigr\|^2_{L^2(\omega)}\le {C\over \delta}{\cal
E}(u_\delta,\Omega_\delta)\le C,\leqno (2.10)$$ hence
$$\Bigl\|{\partial {\cal U}_{3,\delta}\over \partial x_1}+b_{2,\delta}\Bigr\|^2_{L^2(\omega)}+\Bigl\|{\partial {\cal U}_{3,
\delta}\over \partial x_2}-b_{1,\delta}\Bigr\|^2_{L^2(\omega)}\le {C\over \delta^3}{\cal E}(u_\delta,\Omega_\delta)\le
{C\over \delta^2}.$$ Now we apply  the Poincar\'e-Wirtinger inequality to the function ${\cal
U}_{3,\delta}+ b_{2,\delta}x_1-b_{1,\delta}x_2$. There exists $a_{3,\delta}$ such that 
$$\displaystyle ||{\cal U}_{3,\delta}-a_{3,\delta}+ b_{2,\delta}x_1-b_{1,\delta}x_2||^2_{L^2(\omega)}\le  {C\over
\delta^3}{\cal E}(u_\delta,\Omega_\delta)\le{C\over \delta^2}.$$ The sequences ${\cal U}_{M,\delta}-r_{M,\delta}$,
$\delta\bigl\{{\cal U}_{3,\delta}-a_{3,\delta}+b_{2,\delta}x_1- b_{1,\delta}x_2\bigr\}$, $\delta({\cal
R}_{\alpha,\delta}-b_{\alpha,\delta})$ and
$\displaystyle{1\over \delta}{\cal T}_\delta(
\widetilde{u}_\delta)$ are bounded in $H^1(\omega,\R^2)$ (respectively $H^1(\omega)$ and
$L^2(\omega, H^1(]-1,1],\R^3))$. We extract from these sequences some subsequences, still denoted in the same way, such
that
$$\left\{\eqalign{ &{\cal U}_{1,\delta}-a_{1,\delta}+x_2 b_{3,\delta}\rightharpoonup U_1,\quad {\cal
U}_{2,\delta}-a_{2,\delta}-x_1 b_{3,\delta}\rightharpoonup U_2\quad \hbox{weakly in}\quad H^1(\omega),\cr
& \delta\bigl\{{\cal U}_{3,\delta}-a_{3,\delta}+x_1b_{2,\delta}-x_2b_{1,\delta}\bigr\}\rightharpoonup U_3\quad
\hbox{weakly in}\quad H^1(\omega),\cr  
&\delta({\cal R}_{\alpha,\delta}-b_{\alpha,\delta})\rightharpoonup {\cal R}_\alpha\quad \hbox{weakly in}\quad 
H^1(\omega),\cr 
&{1\over \delta}{\cal T}_\delta(\widetilde{u}_\delta) \rightharpoonup \widetilde{u}\qquad\hbox{weakly in}\quad L^2(\omega,
H^1(]-1,1[,\R^3)).\cr}\right.\leqno (2.11)$$ The limits of the sequences
$\displaystyle \delta\Bigl\{ {\partial {\cal U}_{3,\delta}\over \partial x_1}+{\cal R}_{2,\delta}\Bigr\}$ and 
$\displaystyle\delta\Bigl\{{\partial {\cal U}_{3,\delta}\over \partial x_2}-{\cal R}_{1,\delta}\Bigr\}$ are equal to zero  by
$ (2.10)$, hence the equalities 
$${\partial U_3\over \partial x_1}=-{\cal R}_2,\qquad {\partial U_3\over \partial x_2}={\cal R}_1,\leqno (2.12)$$ and
the belonging of $U_3$ to $H^2(\omega)$. From the limits $ (2.11)$ and from the equalities $ (2.12)$ we
immediately  deduce the limits of the unfolded ${\cal T}_\delta\bigl(u_{1,\delta}-a_{1,\delta}+x_2 b_{3,\delta}
\bigr)$, ${\cal T}_\delta\bigl(u_{2,\delta}-a_{2,\delta}-x_1 b_{3,\delta}\bigr)$ and $\delta {\cal
T}_\delta\bigl(u_{3,\delta}-a_{3,\delta}+x_1b_{2,\delta}-x_2b_{1,\delta}\bigr)$ in $ H^1(\Omega)$.

\noindent To calculate the components of the strain tensor we use the equality $ (2.3)$
$$\eqalign{ &\gamma_{\alpha\beta}(u_\delta)={1\over 2}\Bigl\{{\partial{\cal U}_{\alpha,\delta}\over \partial
x_\beta}+{\partial {\cal U}_{\beta,\delta}\over \partial x_\alpha}\Bigr\}-x_3{\partial^2{\cal U}_{3,\delta}\over \partial
x_\alpha\partial x_\beta}+{1\over 2}\Bigl\{{\partial \widetilde{u}_{\alpha,\delta}\over \partial
x_\beta}+{\partial\widetilde{u}_{\beta,\delta}\over \partial x_\alpha}\Bigr\},\cr 
&\gamma_{\alpha3}(u_\delta)={1\over 2}\Bigl\{{\partial\widetilde{u}_{\alpha,\delta}\over \partial x_3}
+{\partial\widetilde{u}_{3,\delta}\over \partial x_\alpha}\Bigr\},\qquad
\gamma_{33}(u_\delta)={\partial \widetilde{u}_{3,\delta}\over \partial x_3}.\cr}$$ These equalities are transformed
 through unfolding. All the sequences ${\cal T}_\delta\bigl(\gamma_{ij}(u_\delta)\bigr)$ are bounded in 
$L^2(\Omega)$ and they have a limit in $H^{-1}(\Omega)$, which can be explained thanks to the  convergences $ (2.11)$ and the
equalities $(2.12)$. Hence the last limits of $ (2.9)$.\fin
\noindent{\bf Remark 2.7 : } We consider again the sequence of displacements  $\bigl(u_\delta\bigr)_{\delta>0}$ of Theorem
2.7. We put $\bigl(U^{'}_{e,\delta}\bigr)_{\delta>0}$ another sequence of e.p.d. verifying 
$${\cal E}(U^{'}_{e,\delta},\Omega_\delta)+{\cal D}(u_\delta-U^{'}_{e,\delta},\Omega_\delta)+{1\over \delta^2}
||u_\delta-U^{'}_{e,\delta}||^2_{L^2(\Omega_\delta, \R^3)}\le C{\cal E}(u_\delta, \Omega_\delta)\le C\delta,\leqno (2.13)$$ where
the constant is independent on $\delta$. Then we obtain
$$ ||{\cal U}_{\delta}-{\cal U}^{'}_{\delta}||_{L^2(\omega, \R^3)}\le C\delta\qquad 
 ||{\cal R}_{\alpha,\delta}-{\cal R}^{'}_{\alpha,\delta}||_{L^2(\omega)}\le C.$$ The displacement $u_\delta$  is decomposed
now into the sum of a new Kirchhoff-Love displacement and a new residual one 
$\widetilde{u}^{'}_\delta$
$$u_\delta=\Bigl({\cal U}^{'}_{1,\delta}-x_3{\partial {\cal U}^{'}_{3,\delta}\over \partial x_1}\Bigr)\vec e_1+
\Bigl({\cal U}^{'}_{2,\delta}-x_3{\partial {\cal U}^{'}_{3,\delta}\over \partial x_2}\Bigr)\vec e_2+ {\cal U}^{'}_{3,\delta}
\vec e_3+\widetilde{u}^{'}_\delta$$ The displacement  $\widetilde{u}^{'}_\delta$ verifies the inequality
$\displaystyle ||\widetilde{u}^{'}_\delta||^2_{L^2(\Omega_\delta,\R^3)}+\delta^2\Bigl\|{\partial
\widetilde{u}^{'}_\delta\over\partial x_3}\Bigr\|^2_{L^2(\Omega_\delta,\R^3)}\le C\delta^3.$ After extraction of
subsequences  expressed by the same notation, we obtain the convergences
$$\left\{\eqalign{ 
&{\cal U}_\delta-{\cal U}^{'}_\delta\longrightarrow 0 \qquad \delta\bigl({\cal R}_\delta-{\cal R}^{'}_\delta\bigr)
\longrightarrow 0\qquad  \hbox{strongly in}\quad L^2(\omega,\R^3)\cr 
&{1\over \delta}{\cal T}_\delta(\widetilde{u}^{'}_\delta)\rightharpoonup \widetilde{u}^{'}\qquad
\hbox{weakly in}\quad L^2(\omega,H^1(]-1,1[,\R^3)).\cr}\right.$$ The limits of the unfolded
${\cal T}_\delta\bigl(\gamma_{i3}(u_\delta)\bigr)$ give
$\displaystyle {\partial\widetilde{u}\over\partial t_3}={\partial \widetilde{u}^{'}\over\partial t_3}$. Except for the limit of the
sequence of the unfolded ${1/ \delta}{\cal T}_\delta(\widetilde{u}_\delta)$, the limits $ (2.8)$ and $ (2.9)$  do not depend on
the decomposition of the displacement $u_\delta$ into the sum of an e.p.d. and a residual displacement. What matters is to be
able to approximate $u_\delta$ with the help of an e.p.d. that verifies the estimates $ (2.13)$. It is to be noticed that the mere
knowing of the limits of the unfolded of the stain tensor components of the sequence $\bigl(u_\delta\bigr)_{\delta>0}$ is not
enough to determine completely the residual displacement
$\widetilde{u}$. It is obtained but for a function of  $L^2(\omega,\R^3)$.
\vskip 1mm
\noindent{\bf 2.3. A second  decomposition of a plate displacement}
\vskip 1mm
\noindent   We consider now a round-rimmed plate $\Omega^{'}_\delta$ with a middle surface  $\omega_\delta$. We denote
$$\eqalign{
\Omega^{'}_\delta&=\bigl\{x\in \R^3\; |\; dist(x,  \omega)<\delta\; \bigr\},\qquad \Gamma_\delta=\bigl\{x\in\R^3\, |\,
dist(x,\partial\omega)<\delta)\bigr\}\cr
\widetilde{\Omega}_\delta&=\bigl\{x\in \R^2\; |\; dist(x,  \omega)<2\delta\; \bigr\}\times
]-\delta,\delta[=\omega_{2\delta}\times ]-\delta,\delta[.\cr}$$
\noindent{\bf Lemma 2.8 : }{\it For any $\delta\in]0,\delta_0]$, there exists an extension operator ${\cal P}_\delta$,
linear and continuous from $H^1(\Omega^{'}_\delta,\R^3)$ into $H^1(\widetilde{\Omega}_\delta,\R^3)$ such that
$${\cal P}_\delta(u)_{|\Omega^{'}_\delta}=u,\qquad 
\left\{\eqalign{ {\cal E}({\cal P}_\delta(u),\widetilde{\Omega}_\delta)\le C{\cal E}(u,\Omega^{'}_\delta)\cr {\cal E}({\cal
P}_\delta(u),\widetilde{\Omega}_\delta\setminus \Omega_\delta)\le C{\cal E}(u,\Gamma_\delta)\cr}\right.
\leqno(2.14)$$  The constants do not depend on  $\delta$.}

\noindent{\bf Proof : } See Annex B.\fin
\noindent  The extension of $u$ to $\widetilde{\Omega}_{\delta}$ is still denoted $u$.
\vskip 1mm
\noindent{\bf A second elementary plate displacement associated with a displacement of  $H^1(\Omega^{'}_\delta, \R^3)$.}
\vskip 1mm
\noindent{\bf Definition 2.9 : } With any $u\in H^1(\Omega^{'}_\delta, \R^3)$  we associate the e.p.d. $U^{'}_e$ defined by 
$$\left\{\eqalign{ 
U^{'}_e(x)&={\cal U}^{'}(\widehat{x})+{\cal R}^{'}(\widehat{x})\land x_3\vec e_3, \qquad x\in \Omega^{'}_\delta,\cr
{\cal U}^{'}(\widehat{x})&={6\over \pi \delta^3}\int_{B(\hat x ; \delta/2)}u(M)dM,\qquad  {\cal R}^{'}(\widehat{x})={60\over
\pi \delta^5}\int_{B(\hat x; \delta/2)}\overrightarrow{\widehat{x}M}\land u(M)dM,\qquad \widehat{ x}\in
\omega_\delta.\cr}\right.\leqno (2.15)$$
\noindent{\bf Theorem 2.10 : }{\it We have the following inequalities:
$$\left\{\eqalign{ &\delta^3||\nabla{\cal R}^{'}||^2_{L^2(\omega_\delta,\R^6)}+\delta\Bigl\|{\partial {\cal U}^{'}\over \partial
x_\alpha}- {\cal R}^{'}\land \vec e_\alpha\Bigr\|^2_{L^2(\omega_\delta,\R^3)}+{\cal E}(U^{'}_e,\Omega^{'}_\delta)\le
C{\cal E}(u,\Omega^{'}_\delta),\cr  &{\cal D}(u-U^{'}_e,\Omega^{'}_\delta)\le C{\cal E}(u,\Omega^{'}_\delta) ,\quad
||u-U^{'}_e||^2_{L^2(\Omega^{'}_\delta, \R^3)}\le C\delta^2{\cal E}(u,\Omega^{'}_\delta).\cr}\right.
\leqno (2.16)$$ The constants  depend only on $\omega$.}  
\vskip 1mm
\noindent{\bf Proof : } We now consider the covering  $\bigl\{\omega^{'}_{\delta,n}\bigr\}_{n\in N_\delta}$  (see 
Lemma 4.1 in Annex A). We put $\omega^{''}_{\delta,n}=\bigl\{\widehat{x}\in
\omega^{'}_{\delta,n}\; |\; dist(\widehat{x},\partial \omega^{'}_{\delta,n})>{\delta/2}\bigr\}$ and  ${\cal O}^{'}_{\delta,n}
=\omega^{'} _{\delta,n}\times ]-\delta,\delta[$, ${\cal O}^{''}_{\delta,n}=\omega^{''} _{\delta,n}\times ]-\delta,\delta[$,
$n\in N_\delta$. The family $\bigl\{\omega^{''}_{\delta,n}\bigr\}_{1\le n\le N_\delta}$ verifies
$$measure\Bigl(\bigcup_{n\in\N_\delta}\omega^{''}_{\delta,n}\setminus \omega_\delta\Bigr)=0$$ From Lemma
2.3 in [4] there exists a rigid displacement $r_n$ such that 
$${\cal D}(u-r_n,{\cal O}^{'}_{\delta,n})+{1\over \delta^2} ||u-r_n||^2_{L^2({\cal O}^{'}_{\delta, n} , \R^3)}\le C{\cal
E}(u,{\cal O}^{'}_{\delta,n}),\quad r_n(x)=a_n+b_n\land\overrightarrow{A_nx}, \quad (a_n,b_n)\in\R^3.\leqno (2.17)$$ 
\noindent The constants do not depend on $n$ nor $\delta$.  We calculate the mean of $(u-r_n)(M)$ and of
$\overrightarrow{\widehat{x}M}\land (u-r_n)(M)$ in the ball $B(\widehat{x},\delta/2)$,
$\widehat{x}\in \omega^{''}_{\delta,n}$. Due to $ (2.17)$ we obtain 
$$\left\{\eqalign{ 
&||{\cal U}^{'}(\widehat{ x})-a_n-b_n\land\overrightarrow{A_n\widehat{x}}||^2_2 \le {C\over \delta} {\cal
E}(u,{\cal O}^{'}_{\delta, n})\;  \Longrightarrow  \;   
||{\cal U}^{'}-a_n-b_n\land\overrightarrow{A_n\cdot}||^2_{L^2(\omega^{''}_{\delta,n},\R^3)} \le C\delta {\cal E}(u,{\cal
O}^{'}_{\delta, n}),\cr 
&\hbox{and}\quad  ||{\cal R}^{'}(\widehat {x})-b_{n}||^2_2\le {C\over\delta^3}{\cal E}(u,{\cal O}^{'}_{\delta,n})\;
\Longrightarrow\ \;  ||{\cal R}^{'}-b_{n}||^2_{L^2(\omega^{''}_{\delta,n},\R^3)}\le {C\over\delta}{\cal E}(u,{\cal
O}^{'}_{\delta,n}),\cr}\right.\leqno (2.18)$$ where $||\cdot||_2$ refers to  the euclidian  norm of $\R^3$.  From the inequalities  $
(2.17)$, $ (2.18)$ and after elimination of the rigid displacement $r_n$,  we obtain   
$ ||u-U^{'}_e||^2_{L^2({\cal O}^{''}_{n,\delta},\R^3)}\le C\delta^2{\cal E}(u,{\cal O}^{'}_{n,\delta}).$ We add all these
inequalities to obtain
$||u-U^{'}_e||^2_{L^2(\Omega^{'}_\delta,\R^3)}\le C\delta^2{\cal E}(u,\Omega^{'}_\delta)$.

\noindent The components ${\cal U}^{'}$ and ${\cal R}^{'}$ of $U^{'}_e$ belong to $H^1(\omega_\delta,\R^3)$. The partial
derivatives of these functions are
$${\partial{\cal U}^{'}\over \partial x_\alpha}(\widehat{x})={6\over \pi \delta^3}\int_{B(\hat x ; \delta/2)}{\partial u\over
\partial x_\alpha}(M)dM,\qquad\enskip {\partial{\cal R}^{'}\over \partial x_\alpha}(\widehat{x})={24\over \pi
\delta^5}\int_{B(\hat x; \delta/2)}\overrightarrow{\widehat{x}M}\land {\partial u\over \partial x_\alpha}(M)dM.$$ Let us
calculate the means of $\displaystyle {\partial \over \partial x_\alpha}(u-r_n)(M)$ and of $\displaystyle
\overrightarrow{\widehat{x}M}\land {\partial \over\partial x_\alpha}(u-r_n)(M)$ in the ball
$B(\widehat{x},\delta/2)$, $\widehat{x}\in \omega^{''}_{\delta,n}$. Thanks to $ (2.17)$ we obtain 
$$\left\{\eqalign{ &\Bigl\|{\partial{\cal U}^{'}\over\partial x_\alpha}(\widehat{x})-b_n\land\vec e_\alpha\Bigr\|^2_2
+\delta^2\Bigl\|{\partial{\cal R}^{'}\over\partial x_\alpha}(\widehat{x})\Bigr\|^2_2\le {C\over \delta^3}{\cal E}(u,{\cal
O}^{'}_{\delta,n})\cr
\Longrightarrow \quad &\Bigr\|{\partial{\cal U}^{'}\over\partial x_\alpha}-b_n\land\vec
e_\alpha\Bigr\|^2_{L^2(\omega^{''}_{\delta,n},\R^3)} +\delta^2\Bigr\|{\partial{\cal R}^{'}\over\partial
x_\alpha}\Bigr\|^2_{L^2(\omega^{''}_{\delta,n},\R^3)}\le {C\over \delta}{\cal E}(u,{\cal O}^{'}_{\delta,n}),\cr}\right.\leqno
(2.19)$$ hence, on the one hand $\displaystyle \delta\Bigr\|{\partial{\cal U}^{'}\over\partial x_\alpha}-{\cal R}^{'}\land\vec
e_\alpha\Bigr\|^2_{L^2(\omega^{''}_{n,\delta},\R^3)}\le C{\cal E}(u,{\cal O}^{'}_{n,\delta})$ using $(2.18)$ and  on the
other hand $\displaystyle  \delta^3\|\nabla{\cal R}\|^2_{L^2(\omega^{''}_{n,\delta},\R^6)}\le  C{\cal E}(u,{\cal
O}^{'}_{n,\delta})$ . We add all these inequalities 
$$\delta\Bigr\|{\partial{\cal U}^{'}\over\partial x_\alpha}-{\cal R}^{'}\land\vec
e_\alpha\Bigr\|^2_{L^2(\omega_\delta,\R^3)}\le C{\cal E}(u,\Omega^{'}_\delta),\qquad\delta^3\|\nabla{\cal
R}\|^2_{L^2(\omega_\delta,\R^6)}\le  C{\cal E}(u,\Omega^{'}_{\delta}).$$ From $ (2.18)$ and $(2.19)$ once more we obtain  
$${\cal D}(u-U^{'}_e,{\cal O}^{''}_{n,\delta})\le C{\cal E}(u,{\cal O}^{'}_{n,\delta})\quad \Longrightarrow\quad {\cal
D}(u-U^{'}_e,\Omega^{'}_{\delta})\le C{\cal E}(u,\Omega^{'}_{\delta}).$$  Theorem 2.10 is proved.\fin
\noindent{\bf 3. The displacements of a structure and the asymptotic behavior of a structure made of plates}
\vskip 1mm
\noindent {\bf 3.1. The structure made of plates} 
\vskip 1mm
\noindent  We work on a set of $N$ plane bounded domains  with polygonal boundary,  included in $\R^3$,
$\bigl(\omega_l\bigr)_{1\le l\le N}$. The skeleton $\GS$ is the union of  $\bigl(\overline{\omega}_l\bigr)_{1\le l\le N}$. A
{\it face} of $\GS$ is a closed set
$\overline{\omega}_l$. An {\it edge} of $\GS$ is a maximal segment shared by a set of faces or a
maximal segment belonging to the boundary of a face. A {\it vertex} of $\GS$ is an extremity of an edge.  

\noindent {\bf Hypotheses : } We suppose 

 {\bf H1} $\bullet$ for any pair of faces $(\overline{\omega}_l,\overline{\omega}_p)$, there exists a sequence of faces 
$\overline{\omega}_l=\overline{\omega}_{l_0},\enskip\overline{\omega}_{l_1},\enskip\ldots$
$\overline{\omega}_{l_k} =\overline{\omega}_p$ such that $\overline{\omega}_{l_r}$ and $\overline{ \omega}_{l_{r+1}}$
have an edge in common, $0\le r\le k-1$,

{\bf H2} $\bullet$  for any vertex $A$ and any pair of faces $(\overline{\omega}_l,\overline{\omega}_p)$ containing
$A$, there exists a sequence of faces 
$\overline{\omega}_l=\overline{\omega}_{l_0},\enskip\overline{\omega}_{l_1},\enskip\ldots
\enskip\overline{\omega}_{l_k} =\overline{\omega}_p$ such that $\overline{\omega}_{l_r}$ and $\overline{ \omega}_{l_{r+1}}$
have an edge in common containing $A$, $0\le r\le k-1$,

{\bf H3} $\bullet$   the skeleton $\GS$  is fixed all along some  edges.

$$\qquad\vbox{\centerline{\dessin 162.8mm by 81.43mm (Fig1 scaled 600)}} $$
\centerline{Figure 1.  The skeleton $\GS$.}
\vskip 1mm
\noindent We denote

$\bullet$  $\Gamma_0$  the fixed part of the skeleton,

$\bullet$  ${\cal J}$ the set of edges common to several faces,

$\bullet$ ${\cal N}$ the set of vertexes common to several faces.
\vskip 1mm
\noindent  The structure made of plates is the domain ${\cal S}_\delta=\bigl\{x\in \R^3\; |\; dist(x, \GS)<\delta\bigr\}$.  This
structure is made of the  gathering of the plates $\Omega^{'}_{l,\delta}$ with  thickness   $2\delta$, with  middle surface 
$\omega_{l,\delta}$ and with rounded rim. Each domain  $\Omega^{'}_{l,\delta}$ is equipped with a local frame $(O^{(l)};\vec
e^{(l)}_1,\vec e^{(l)}_2,\vec e^{(l)}_3)$, $O^{(l)}\in \omega_l$, $\vec e^{(l)}_3$ is the normal direction  to  the face
 $\overline{\omega}_{l}$. The plate $\Omega^{'}_{l,\delta}$ contains the plate $\Omega_{l,\delta}=\omega_l\times
]-\delta,\delta[$ ($\R^3$ being equipped with the above local frame). The reference plate is the open set
$\Omega_l=\omega_l\times ]-1,1[$ obtained through the transformation of $\Omega_{l,\delta}$ by the orthogonal affiniy
of ratio  ${1/\delta}$.

\noindent 
\noindent The structure ${\cal S}_\delta$ is fixed to part $\Gamma_{0,\delta}=\bigl\{x\in \partial{\cal S}_\delta\; |\;
dist(x,\Gamma_0) =\delta\bigr\}$ of its boundary . For each edge $J\in{\cal J}$ we choose a unit vector  $\vec e_J$ in the
 direction of the edge.

\noindent We consider $\GS_\delta=\displaystyle\bigcup_{l=1}^N\overline{\omega}_{l,\delta}$,  the set of the middle
surfaces of the plates.

\noindent There exist two constants $\delta_0>0$ and $\mu_0>0$  depending only on the skeleton $\GS$ such that for any
$\delta\in ]0,\delta_0]$ the parts common to several plates are in the union of the junctions 
$$\bigcup_{J\in{\cal J}}\Bigl\{x\in \R^3\; |\; dist(x,J)< \eta_0\delta\Bigr\}.$$  
\noindent The  restriction of a function $\phi$, defined on $\GS$, (resp. $\GS_\delta$, ${\cal S}_\delta$), to $\omega_l$, (resp.
$\omega_{l,\delta}$, $\Omega^{'}_{l,\delta}$) is denoted $\phi^{(l)}$. In the same way, we denote 
$x^{(l)}_\alpha=x^{(l)}\cdot\vec e^{(l)}_\alpha$ and $x^{(l)}_3=x^{(l)}\cdot\vec e^{(l)}_3$ the local variables    
$\bigl(x^{(l)}=(\widehat{x}^{(l)},x^{(l)}_3)=(x^{(l)}_1, x^{(l)}_2, x^{(l)}_3)\bigr)$.

\noindent The  space $H^1(\GS,\R^n)$ (resp. $H^1(\GS_\delta,\R^n)$) is the set of the functions $\phi$ defined a.e. in
$\GS$ (resp. $\GS_\delta$), with values in $\R^n$, such that the restriction $\phi^{(l)}$ belongs to $H^1(\omega_l,\R^n)$
(resp. $H^1(\omega_{l,\delta},\R^n)$) and such that for any edge $J\in {\cal J}$ and any pair of faces
$(\overline{\omega}_l,\overline{\omega}_k)$ containing $J$, we have the equality of  the restrictions to J,
$(\phi^{(l)})_{|J}=(\phi^{(k)})_{|J}$ in $H^{1/2}(J,\R^n)$. 
\vskip 1mm
\noindent The space $H^1_{\Gamma_0}(\GS,\R^n)$ is the subspace of $H^1(\GS,\R^n)$ the elements of which are a.e. equal
to zero on  $\Gamma_0$. \noindent We equip $H^1(\GS_\delta,\R^3)$  with the inner product
$$\displaystyle[{\cal U},{\cal V}]=\sum_{l=1}^N\int_{\omega_{l,\delta}}\Bigl\{{\cal U}^{(l)}\cdot {\cal V}^{(l)}
+{\partial {\cal U}^{(l)}\over\partial x^{(l)}_\alpha}\cdot{\partial {\cal V}^{(l)}\over\partial x^{(l)}_\alpha}\Bigr\}.$$  The
associated norm is denoteted $||\cdot||$. 
\vskip 1mm
\noindent{\bf 3.2. The elementary displacements of plates structure} 
\vskip 1mm
\noindent{\bf Definition 3.1 : } An {\bf elementary displacement of a plates structure } (e.d.p.s.) is a
displacement
$\Phi\in H^1({\cal S}_\delta,\R^3)$, such that there exist two elements ${\cal A}$ and ${\cal B}$
belonging to
$H^1(\GS_\delta,\R^3)$,  such that for any $l\in\{1,\ldots, N\}$, 
$$\Phi^{(l)}(x)={\cal A}^{(l)}(\widehat{x}^{(l)})+{\cal B}^{l)}(\widehat{x}^{(l)})\land x^{(l)}_3\vec e^{(l)}_3$$ is an e.p.d. 
of the plate
$\Omega^{'}_{l,\delta}$.
\vskip 1mm
\noindent The functions  ${\cal A}$ and ${\cal B}$ are respectively  the first component and the second component
of the e.d.p.s. $\Phi$. The function ${\cal A}$ accounts for the displacement  of the skeleton faces while 
${\cal B}$ accounts for the rotation of the normal directions to the plates and of the rotation of the
faces around the edges. 
\medskip
\noindent{\bf Theorem 3.2  : }{\it For any displacement $u\in H^1_{\Gamma_{0}}({\cal S}_\delta,\R^3)$ there exists an
e.d.p.s. $U_e$ of components
$({\cal U},{\cal R})\in H^1_{\Gamma_0}(\GS_\delta,\R^3)\times H^1_{\Gamma_0}(\GS_\delta,\R^3) $ such that
$$\left\{\eqalign{ &\sum_{l=1}^N\Bigl\{\delta^3||\nabla{\cal
R}^{(l)}||^2_{L^2(\omega_{l,\delta},\R^6)}+\delta\Bigl\|{\partial {\cal U}^{(l)}\over
\partial x^{(l)}_\alpha}-{\cal R}^{(l)}\land\vec e^{(l)}_\alpha\Bigr\|^2_{L^2(\omega_{l,\delta},\R^3)}\Bigr\} \le C{\cal
E}(u,{\cal S}_\delta),\cr  &{\cal E}(U_e,{\cal S}_\delta)+{\cal D}(u-U_e,{\cal S}_\delta)\le C{\cal E}(u,{\cal S}_\delta),\quad
||u-U_e||^2_{L^2({\cal S}_\delta, \R^3)}\le C\delta^2{\cal E}(u,{\cal S}_\delta).\cr}\right.\leqno (3.1)$$}
\noindent{\bf Proof : } Let $u$ be in $H^1_{\Gamma_{0}}({\cal S}_\delta,\R^3)$. Thanks to Lemma 4.1, for any $l$ belonging
to $\{1,\ldots,N\}$, we extend  the restriction $u^{(l)}$  to the plate  $\Omega^{'}_{l,\delta}$, into a displacement  of
$\widetilde{\Omega}_{l,\delta}$  (the plate of thickness $2\delta$  and of middle surface $\omega_{l,2\delta}$). Therefore,
using the formulas  $ (2.15)$, we can define an   e.p.d. $U_e^{'(l)}$ of the plate $\Omega^{'}_{l,\delta}$ verifying $ (2.16)$.
Both components ${\cal U}^{'(l)}$ and
${\cal R}^{'(l)}$ of $U_e^{'(l)}$ are the restrictions  of elements belonging to 
$H^1_{\Gamma_0}(\GS_\delta,\R^3)$ (corollary of Lemma 2.8 in Annex B).
\vskip 1mm
 \noindent Then we build a new e.d.p.s. $U_e$ equal to $U_e^{'}$ in the open set $\displaystyle
\bigcup_{J\in {\cal J}}\bigl\{x\in{\cal S}_\delta\; |\; dist(x,J)>2\eta_0\delta\bigr\}$ and equal to an elementary displacement
of rods structures in the  junctions (see Annex B). We then can deduce $ (3.1)$.\fin
\noindent{\bf Proposition 3.3 (Korn inequality) : }{\it For any displacement $u\in H^1_{\Gamma_{0}}({\cal S}_\delta,\R^3)$, we
have
$$\delta||{\cal R}||^2+\delta||{\cal U}||^2+{\cal D}(u,{\cal S}_\delta)+||u||^2_{L^2({\cal S}_\delta,\R^3)}\le
{C\over \delta^2}{\cal E}(u,{\cal S}_\delta). \leqno (3.2)$$ The constant does not depend on $\delta$.}

\noindent{\bf Proof : } The estimates $ (3.1)$ of the gradients  of the  ${\cal R}^{(l)}$  functions, the nullity of
${\cal R}$ on $\Gamma_0$ and the  hypothesis {\bf H1}   allow us to obtain step by step $||{\cal
R}||^2_{L^2(\GS_\delta,\R^3)}\le C/\delta^3{\cal E}(u,{\cal S}_\delta)$. This inequality and $ (3.1)$ give then an upperbound
of the  $L^2$ norms of the functions gradients  ${\cal U}^{(l)}$. The nullity of ${\cal U}$ on
$\Gamma_0$ and the hypothesis {\bf H1} imply then that $||{\cal U}||^2_{L^2(\GS_\delta,\R^3)}\le C/\delta^3{\cal
E}(u,{\cal S}_\delta)$. From these estimates of ${\cal U}$ and ${\cal R}$ follow   ${\cal D}(U_e,{\cal
S}_\delta)\le {C/ \delta^2}{\cal E}(u,{\cal S}_\delta)$ and $||U_e||^2_{L^2({\cal S}_\delta,\R^3)}\le {C/ \delta^2}{\cal E}(u,{\cal
S}_\delta)$. Then again, thanks to  $(3.1)$,  we obtain the estimates of the $L^2$ norm  of $u$ and of its gradient.\fin
\noindent{\bf 3.3. Inextensional displacements, extensional displacements}  
\vskip 1mm 
\noindent The space $H^1_{\rho,\Gamma_0}(\GS,\R^3)$  is the set of the functions $\phi$ defined a.e. in $\GS$, with values
in $\R^3$, such that:

$\bullet$  for any $l\in\{1,\ldots, N\}$, the restrictions $\phi^{(l)}_1$  and $\phi^{(l)}_2$ belong to
$H^1(\omega_{l})$ and   $\phi^{(l)}_3$ belongs to
$$H^1_\rho(\omega_{l})=\Bigl\{\psi\in L^2(\omega_{l})\; |\; \sqrt{\rho}\,\nabla\psi\in [L^2(\omega_{l})]^2\;\Bigr\}$$
where $\rho(\widehat{x})=dist(\widehat{x},{\cal N})$ (distance from the point  $\widehat{x}\in \GS_\delta$ to the vertexes
belonging to several faces),  

$\bullet$ for any edge $J\in {\cal J}$ and any pair of faces $(\overline{\omega}_l,\overline{\omega}_k)$ containing
$J$, we have $(\phi^{(l)})_{|J}=(\phi^{(k)})_{|J}$ in  $H^{1/2}(J,\R^3)$,

$\bullet$ the function $\phi$ is equal to zero on $\Gamma_0$.
\vskip 1mm
\noindent  We equip  $H^1_{\rho,\Gamma_0}(\GS,\R^3)$ with the inner product
$$<{\cal U},{\cal V}>_{\rho}=\sum_{l=1}^N\int_{\omega_{l}}\Bigl\{\gamma_{\alpha\beta}({\cal U}^{(l)})\gamma_{
\alpha\beta}({\cal V}^{(l)}) +\rho\nabla {\cal U}^{(l)}_3\cdot\nabla {\cal V}^{(l)}_3\Bigr\},$$  and with the norm $|{\cal
U}|_{\rho}=\sqrt{<{\cal U},{\cal U}>_{\rho}}$. The usual norm on $H^1_{\rho,\Gamma_0}(\GS,\R^3)$ is
$$||{\cal U}||_{\rho}=\sqrt{\sum_{l=1}^N\int_{\omega_{l}}\Bigl\{|\nabla{\cal U}^{(l)}_1|^2  +|\nabla{\cal U}^{(l)}_2|^2
+\rho |\nabla {\cal U}^{(l)}_3|^2\Bigr\}}$$ 
\noindent{\bf Lemma 3.4  : }{\it The norms $||\cdot||_{\rho}$ and $|\cdot|_{\rho}$ are equivalent in $
H^1_{\rho,\Gamma_0}(\GS,\R^3)$.  Moreover the space 
$H^1_{\Gamma_0}(\GS,\R^3)$ is dense in  $H^1_{\rho,\Gamma_0}(\GS,\R^3)$.}

\noindent{\bf Proof : } See Annex C.\fin
\noindent{\bf Definition 3.5 : } An {\bf inextensional displacement of the skeleton} is an element $U$ belonging to
$H^1_{\rho,\Gamma_0}(\GS,\R^3)$ such that
$$ \forall l\in\{1,\ldots,N\},\qquad \gamma_{\alpha\beta}(U^{(l)})=0,\qquad \hbox{in}\qquad \omega_{l}.$$
\noindent The membrane component $U^{(l)}_M=U^{(l)}_1\vec e^{(l)}_1+U^{(l)}_2\vec e^{(l)}_2$ of an inextensional
displacement is a  rigid displacement of the face $\overline{\omega}_{l}$. The inextensional displacements space
of the skeleton is  denoted $D_I(\GS)$.
\vskip 1mm
\noindent{\bf Definition 3.6 : } An {\bf extensional displacement of the skeleton} is an element of the orthogonal 
$D_E(\GS)$ of $D_I(\GS)$ in $H^1_{\rho,\Gamma_0}(\GS,\R^3)$.  

\noindent The set of  extensional displacements is equipped with the semi-norm
$$||U||_{E}=\sqrt{\sum_{l=1}^N\int_{\omega_{l}}\gamma_{\alpha\beta}(U^{(l)})\gamma_{\alpha\beta}(U^{(l)})},\qquad U\in
D_E(\GS).$$ The semi-norm $||\cdot||_{E}$ is a norm, because if $U\in D_E(\GS)$ is such that  $||U||_{E}=0$ then,
$\gamma_{\alpha\beta}(U^{(l)})=0$ for any $l$. The displacement $U$ is then of  inextensional type  and is  equal to zero. 
\vskip 1mm
\noindent{\bf Lemma 3.7 : }{\it The norms  $||\cdot||_{E}$ and $|\cdot|_{\rho}$ are equivalent in $D_E(\GS)$.}
\vskip 1mm
\noindent{\bf Proof : } See Annex C.\fin 
\noindent{\bf Corollary of Lemma 3.7  : }  Let $u$ be a displacement belonging to $H^1_{\Gamma_0}({\cal S}_\delta,\R^3)$
and $U_e$ the e.d.p.s. given by Theorem 3.2. The restriction to $\GS$ of the first  component ${\cal U}$ of  $U_e$
can be written as   the sum of an extensional displacement and an inextensional displacement,  
$${\cal U}=U_E+U_I,\qquad  U_E\in D_{E}(\GS),\quad U_I\in D_{I}(\GS).\leqno (3.3)$$  According to the inequalities
$ (3.1)$  and Lemma 3.7 we have
$$|U_E|^2_{\rho}\le C ||U_E||^2_{E}\le {C\over\delta}{\cal E}(u,{\cal S}_\delta),\qquad\quad |U_I|_{\rho}^2\le
{C\over\delta^3}{\cal E}(u,{\cal S}_\delta).\leqno (3.4)$$ The constants are independent of
$\delta$.
\vskip 1mm
\noindent{\bf 3.4. The limit displacements} 
\vskip 1mm 
\noindent Let  $\bigl(u_\delta\bigr)_{\delta>0}$ be a sequence of displacements belonging to
$H^1_{\Gamma_0}({\cal S}_\delta,\R^3)$ and verifying   
$${\cal E}(u_\delta,{\cal S}_\delta)\le C\delta, \leqno(3.5)$$ where the constant is independent of $\delta$. Thanks
to the estimates $(3.2)$, $(3.4)$ and $(2.4)$,  from the sequences $\delta{\cal U}_\delta$, $\delta{\cal R}_\delta$,
$\delta U_{I,\delta}$,  $U_{E,\delta}$ and $\displaystyle{1\over \delta}{\cal T}_\delta(\widetilde{u}^{(l)}_\delta)$ we extract
some sub-sequences, still denoted in the same way and which weakly converge, 
$$\left\{\eqalign{
\delta {\cal U}_\delta&\rightharpoonup U_I,\qquad \delta{\cal R}_\delta\rightharpoonup {\cal R}\qquad\hbox{weakly
in}\quad H^1_{\Gamma_0}(\GS,\R^3),\cr 
\delta U_{I,\delta}&\rightharpoonup U_I\qquad\hbox{weakly in}\quad D_I(\GS),\cr 
U_{E,\delta}&\rightharpoonup U_E\qquad\hbox{weakly in}\quad D_E(\GS),\cr  {1\over \delta}{\cal
T}_\delta(\widetilde{u}^{(l)}_\delta) &\rightharpoonup
\widetilde{u}^{(l)}\qquad\hbox{weakly in}\quad L^2(\omega_l, H^1(]-1,1[,\R^3)).\cr}\right.\leqno(3.6)$$
The sequences $\delta {\cal U}_\delta$ and $\delta U_{I,\delta}$ have the same limit in $H^1_{\rho,\Gamma_0}(\GS,\R^3)$
because the sequence  $\delta U_{E,\delta}$ converges to  0 in
$H^1_\rho(\GS,\R^3)$.  After passing to the limit and from  $ (3.1)$ comes 
$$\forall l\in\{1,\ldots,N\},\quad\forall\alpha\in\{1,2\}, \qquad {\partial U^{(l)}_I\over \partial x^{(l)}_\alpha}={\cal
R}^{(l)}\land\vec e_\alpha^{(l)}.\leqno(3.7)$$ Now we define the space of the inextensional displacements limits. We put
$${\cal D}_I(\GS)=\Bigl\{{\cal A}\in D_I(\GS)\cap H^1_{\Gamma_0}(\GS,\R^3)\; |\; \exists {\cal B}\in
H^1_{\Gamma_0}(\GS,\R^3),\; 
 {\partial {\cal A}^{(l)} \over\partial x_\alpha^{(l)}}= {\cal B}^{(l)}\land \vec e^{(l)}_\alpha,\quad \forall l\in\{1,\ldots,N\}
\Bigr\}$$ For any ${\cal A}\in {\cal D}_I(\GS)$, there is only one ${\cal B}$ which  we denote   $\widehat{\nabla}{\cal
A}$. Then we have
$$\forall l\in\{1,\ldots,N\},\qquad \forall \alpha\in\{1,2\},\qquad {\partial {\cal A}^{(l)} \over\partial
x_\alpha^{(l)}}=\widehat{\nabla}{\cal A}^{(l)}\land \vec e^{(l)}_\alpha.$$ We equip
${\cal D}_I(\GS)$ with the norm  $||{\cal A}||_I=||\widehat{\nabla}{\cal A}||_{H^1(\GS,\R^3)}$. The inextensional limit
displacement $U_I$ belongs to ${\cal D}_I(\GS)$. 
\vskip 1mm
\noindent{\bf 3.5. Limit of the unfolded  displacements and limit of the  unfolded   strain tensor components}
\vskip 1mm
\noindent Using  $ (2.8)$ and after transformation by unfolding we have the following limits in the refering plates :
$$\left\{\eqalign{ &\delta{\cal T}_\delta(u_\delta^{(l)}) \rightharpoonup U^{(l)}_I\qquad\hbox{weakly in}\quad
H^1(\Omega_l,\R^3)\cr   
&{\cal T}_\delta(u^{(l)}_{\alpha,\delta}-U^{(l)}_{I,\alpha,\delta}) \rightharpoonup
U^{(l)}_{E,\alpha}-t^{(l)}_3{\partial U^{(l)}_{I,3}\over
\partial x^{(l)}_\alpha}\qquad\hbox{weakly in}\quad H^1(\Omega_l),\cr &{\cal
T}_\delta\bigl(\gamma_{\alpha\beta}(u^{(l)}_\delta)\bigr)\rightharpoonup {1\over 2}\Bigl\{{\partial U^{(l)}_{E,
\alpha}\over \partial x^{(l)}_\beta}+{\partial U^{(l)}_{E,\beta}\over \partial x^{(l)}_\alpha}\Bigr\}
-t^{(l)}_3{\partial^2U^{(l)}_{I,3}\over
\partial x^{(l)}_\alpha\partial x^{(l)}_\beta}\quad \hbox{weakly in}\quad L^2(\Omega_l),\cr &{\cal
T}_\delta\bigl(\gamma_{\alpha3}(u^{(l)}_\delta)\bigr)\rightharpoonup {1\over 2}{\partial\widetilde{u}^{(l)}_\alpha
\over\partial t^{(l)}_3},\qquad{\cal T}_\delta\bigl(\gamma_{33}(u^{(l)}_\delta)\bigr)\rightharpoonup {\partial
\widetilde{u}^{(l)}_3\over \partial t^{(l)}_3}\qquad \hbox{weakly in}\quad L^2(\Omega_l).\cr} \right.\leqno(3.8)$$ 
\noindent{\bf 3.6. Elasticity problem}
\vskip 1mm
\noindent The plates are made of an homogeneous and isotropic  material. Our equations are given within the framework of
linearised elasticity. In ${\cal S}_\delta$ let the elasticity system be 
$$\left\{\eqalign{ - {\partial \over \partial x_j}\{\ a_{iji^{'}j^{'}}{\partial u_{i^{'},\delta}\over \partial x_j^{'}}\} &
=F_{i,\delta}\enskip\hbox { in }\enskip{\cal S}_\delta, \cr 
 u_\delta & =0 \hskip 1.4em \hbox { on }\hskip 0.7em\Gamma_{0,\delta},\cr a_{iji^{'}j^{'}} {\partial u_{i^{'},\delta}
\over\partial x_j^{'}}n_j &=0\enskip\hbox { in }\hskip 0.7em 
\Gamma_\delta,\qquad \Gamma_\delta=\partial {\cal S}_\delta
\setminus\Gamma_{0,\delta}.\cr}\right.\leqno(3.9)$$ \noindent The variational formulation of the problem $(3.9)$ is
$$\left\{\eqalign{
 u_\delta & \in H^1_{\Gamma_0}({\cal S}_\delta ,\R^3) \cr 
\int_{{\cal S}_\delta} &a_{ij^{'}j^{'}}\gamma_{ij}(u_\delta)\gamma_{i^{'}j^{'}}(v)=\int_{{\cal S}_\delta}F_\delta\cdot
v\qquad \qquad   \forall v   \in H^1_{\Gamma_0}({\cal S}_\delta;\R^3)\cr} \right.\leqno(3.10)$$ where
$a_{iji^{'}j^{'}}=\lambda\delta_{ij}\delta_{i^{'}j^{'}}+
\mu(\delta_{ii^{'}}\delta_{jj^{'}}+\delta_{ij^{'}}\delta_{ji^{'}}).$ The constants  $\lambda$ and $\mu$ are
the Lam\'e constants of the material. The plates
$\Omega^{'}_{l,\delta}$  are submitted to volume applied forces .  Among these forces we make a
distinction between those  concerning  the extensional displacements and those concerning the inextensional displacements. 
$$F_\delta(x)=\sum_{l=1}^N\bigl\{\delta
f_I(\widehat{x}^{(l)})+f_E(\widehat{x}^{(l)})\bigr\}\hbox{\bf 1}_{\Omega_{l,\delta}}(x),
\qquad f_I,\; f_E\in L^2(\GS,\R^3),\leqno(3.11)$$ where $\Omega_{l,\delta}=\omega_l\times ]-\delta,\delta[$ (in the local
frame ) and where $\hbox{\bf 1}_{\Omega_{l,\delta}}$ is the characteristic function of the open set 
$\Omega_{l,\delta}$. Hence several volume forces are stacked up in the junctions. 

\noindent The function $f_E$ verifies the condition of  orthogonality
$$ \forall V\in D_I(\GS),\qquad \int_\GS f_E\cdot V=0.\leqno(3.12)$$ 
\noindent Let $({\cal U}_\delta,{\cal R}_\delta)$ be the two components of the  e.d.p.s.  associated to the solution
$u_\delta$ of the problem $(3.10)$. In the plate
$\Omega_{l,\delta}$, the displacement $u^{(l)}_\delta$ is the sum of the e.p.d. ${\cal U}^{(l)}_\delta(\widehat{x}^{(l)})+{\cal
R}^{(l)}(\widehat{x}^{(l)})\land x_3^{(l)}\vec e^{(l)}_3$ and a residual displacement. The displacement ${\cal
U}_\delta$ is the sum of an extensional displacement $U_{E,\delta}$ and of an inextensional displacement 
$U_{I,\delta}$. Then, thanks to $ (3.1)$ and $ (3.4)$, we have  
$$\Bigl|{1\over 2\delta}\int_{{\cal S}_\delta}F_\delta\cdot u_\delta-\int_{\GS}f_E\cdot U_{E,\delta}-\delta\int_{\GS}f_I\cdot 
U_{I,\delta}\Bigr|\le C\bigl\{||f_E||_{L^2(\GS,\R^3)}+||f_I||_{L^2(\GS,\R^3)}\bigr\}\sqrt{{\cal E}(u_\delta,{\cal
S}_\delta)}\leqno(3.13)$$ hence 
$$\bigl|\int_{{\cal S}_\delta}F_\delta\cdot u_\delta\bigr|\le C\bigl\{||f_E||_{L^2(\GS,\R^3)}+||f_I||_{L^2(\GS,\R^3)}\bigr\}\sqrt
\delta\sqrt{{\cal E}(u_\delta,{\cal S}_\delta)}.$$ We deduce that the solution of the variational problem
$(3.10)$ verifies the estimation
$${\cal E}(u_\delta,{\cal S}_\delta)\le C\delta\bigl\{||f_E||^2_{L^2(\GS,\R^3)}+||f_I||^2_{L^2(\GS,\R^3)}\bigr\}.\leqno(3.14)$$
\noindent{\bf 3.7. Asymptotic behavior of the stress tensor} 
\vskip 1mm
\noindent We begin with determining the partial derivatives of the residual displacements
$\widetilde{u}^{(l)}$  in the normal directions to the plates.

\noindent Let $\phi$ be a displacement of $H^1(\Omega_l,\R^3)$, equal to zero in the neighborhood of all the sets
$J\times]-1,1[$ where $J$ is an edge of the face $\overline{\omega}_l$.  For $\delta$ small, the displacement
$\displaystyle \phi_\delta(x)=\delta\phi\bigl(\widehat{x}^{(l)},{x^{(l)}_3\over \delta}\bigr)$ is an acceptable displacement of
the full structure  ${\cal S}_\delta$. We have the following strong convergences of the unfolded of the strained tensor components
of
 $\phi_\delta$:
$$\left\{\eqalign{ 
{\cal T}_\delta\bigl(\gamma_{\alpha\beta}(\phi_\delta)\bigr)&\longrightarrow 0 \qquad \hbox{strongly
in}\quad  L^2(\Omega_l),\cr  
{\cal T}_\delta\bigl(\gamma_{\alpha 3}(\phi_\delta)\bigr)&\longrightarrow {1\over
2}{\partial\phi_\alpha \over \partial t^{(l)}_3}\qquad \hbox{strongly in}\quad L^2(\Omega_l),\cr {\cal
T}_\delta\bigl(\gamma_{33}(\phi_\delta)\bigr)&\longrightarrow{\partial\phi_3\over \partial t^{(l)}_3}\qquad
\hbox{srtongly in}\quad L^2(\Omega_l).\cr}\right.\leqno(3.15)$$ 
\noindent  We now take $\phi_\delta$ as a test-displacement in $(3.10)$, we  transform by unfolding the integral on
$\Omega_{l,\delta}$ into an integral on $\Omega_l$ and after dividing by the thickness of the plate we pass to the limit. We
obtain
$$ \int_{\Omega_l}\Bigl[\lambda\Bigl\{{\partial U^{(l)}_{E,1}\over\partial x^{(l)}_1}-t_3^{(l)}{\partial^2 U^{(l)}_{I,3}\over
\partial x^{(l),2}_1}+{\partial U^{(l)}_{E,2}\over\partial x^{(l)}_2}-t_3^{(l)}{\partial^2U^{(l)}_{I,3}\over\partial 
x^{(l),2}_2}\Bigr\}+(\lambda+2\mu){\partial\widetilde{u}^{(l)}_{3}\over\partial t_3^{(l)}}\Bigr]{\partial\phi_3
\over\partial t_3^{(l)}}  +\mu \Bigl[{\partial\widetilde{u}^{(l)}_{1}\over\partial
t_3^{(l)}}{\partial\phi_1\over\partial t_3^{(l)}}+ {\partial\widetilde{u}^{(l)}_{2}\over\partial
t_3^{(l)}}{\partial\phi_2\over\partial t_3^{(l)}}\Bigr]=0\leqno(3.16)$$ because the right member of
$(3.10)$ tends to 0 $\bigl(\displaystyle {1\over 2\delta}\bigl|\int_{\Omega_{l,\delta}}F_\delta\cdot
\phi_\delta\bigr|\le C\delta||\phi||_{L^2(\Omega_l,\R^3)}\bigr)$.

\noindent The set of  these test-displacements is a dense subset in 
$L^2(\omega_l , H^1(]-1,1[,\R^3))$. Hence the equality $(3.16)$ is  verified for any element of  $L^2(\omega_l ,
H^1(]-1,1[,\R^3))$. We deduce the partial derivatives $\displaystyle{\partial\widetilde{u}^{(l)}\over
\partial t^{(l)}_3}$ in terms of the first partial  derivatives of  
$U^{(l)}_E$  and of the second partial derivatives of $U^{(l)}_I$,
$${\partial\widetilde{u}^{(l)}_1\over \partial t_3^{(l)}}={\partial\widetilde{u}^{(l)}_2\over \partial t_3^{(l)}}=0,\qquad
{\partial\widetilde{u}^{(l)}_3 \over
\partial t_3^{(l)}}={\lambda\over \lambda+2\mu}\Bigl\{-{\partial U^{(l)}_{E,1}\over\partial x^{(l)}_1}-{\partial
U^{(l)}_{E,2}\over\partial x^{(l)}_2}+t^{(l)}_3\Delta U^{(l)}_{I,3}\Bigr\}.\leqno(3.17)$$

\noindent We give now the weak limit in $L^2(\Omega_l)$ of the unfolded of the stress tensor components
$$\left\{\eqalign{ {\cal T}_\delta\bigl(\sigma_{11}(u^{(l)}_\delta)\bigr)&\rightharpoonup{E\over 1-\nu^2}\Bigl[{\partial
U^{(l)}_{E,1}\over\partial x^{(l)}_1}-t^{(l)}_3{\partial^2U^{(l)}_{I,3}\over\partial  x^{(l),2}_1}+\nu\Bigl\{{\partial
U^{(l)}_{E,2}\over\partial x_2^{(l)}}-t^{(l)}_3{\partial^2 U^{(l)}_{I,3}\over \partial x^{(l),2}_2}\Bigr\}\Bigr],\cr {\cal
T}_\delta\bigl(\sigma_{12}(u^{(l)}_\delta)\bigr)&\rightharpoonup\mu\Bigl\{{\partial U^{(l)}_{E,1}\over\partial
x_2^{(l)}}+{\partial  U^{(l)}_{E,2}\over\partial x_2^{(l)}}-2 t^{(l)}_3{\partial^2 U^{(l)}_{I,3}\over \partial x^{(l)}_1\partial
x_2^{(l)}}\Bigr\},\cr  {\cal T}_\delta\bigl(\sigma_{22}(u^{(l)}_\delta)\bigr)&\rightharpoonup{E\over 1-\nu^2}\Bigl[{\partial
U^{(l)}_{E,2}
\over\partial x_2^{(l)}}-t^{(l)}_3{\partial^2U^{(l)}_{I,3}\over\partial x_2^{(l),2}}+\nu\Bigl\{{\partial
U^{(l)}_{E,1}\over\partial x_1^{(l)}}-t_3^{(l)}{\partial^2 U^{(l)}_{I,3}\over \partial x^{(l),2}_1}\Bigr\}\Bigr],\cr  {\cal
T}_\delta\bigl(\sigma_{i3}(u^{(l)}_\delta)\bigr)&\rightharpoonup 0.\cr}\right.\leqno(3.18)$$
\noindent{\bf 3.8. The extensional displacement  $U_E$ or the problem of coupled membrane plates}
\vskip 1mm
\noindent{\bf Theorem 3.8 : }{\it   The extensional  displacement $U_E$ is the solution of the variational problem
$${E\over 1-\nu^2}\sum_{l=1}^N\int_{\omega_l}\bigl[(1-\nu) \gamma_{\alpha\beta}(U^{(l)}_{E})\gamma_{\alpha\beta}( 
V^{(l)})+\nu\gamma_{\alpha\alpha}(U^{(l)}_{E})\gamma_{\beta\beta}( V^{(l)})\bigr]=\int_\GS f_E\cdot V,\qquad\forall V\in
D_E(\GS).\leqno(3.19)$$ where $E$ is the Young  modulus and $\nu$ the Poisson constant .\fin}   
\smallskip
\noindent The proof of Theorem 3.8 requires the next lemma.
\smallskip
\noindent{\bf Lemma 3.9 : }{\it For any element $V\in H^1_{\Gamma_0}(\GS_{\delta_0},\R^3)$, there exists a sequence of
displacements $\bigl(V_\delta)_{0<\delta\le \delta_0}$ belonging to
$H^1_{\Gamma_0}(\GS_{\delta_0},\R^3)\cap H^1_{\Gamma_0}({\cal S}_{\delta},\R^3)$ such that
$$V_\delta\longrightarrow V\qquad \hbox{strongly in}\quad H^1_{\Gamma_0}(\GS_{\delta_0},\R^3).\leqno(3.20)$$}
\noindent{\bf Proof  : } See Annex E.\fin 
\noindent{\bf Proof of Theorem 3.8 : } Let $V$ be an element of $H^1_{\Gamma_0}(\GS,\R^3)$, we extend $V$ into an
element, still denoted $V$, of the space $H^1_{\Gamma_0}(\GS_{\delta_0},\R^3)$.  We take  $V_\delta$  as a
test-displacement in $(3.10)$, we transform,  by unfolding, the integrals on the plates  into integrals on the reference plates  then
we divide by $2\delta$. Thus we are led to take into account again and again the neighborhoods of the edges belonging to
${\cal J}$. 

\noindent Let $J$ be an edge common to several faces. For any face $
\overline{\omega}_l$ containing
$J$, we have
$$\hbox{\bf 1}_{\bigl\{\hat x^{(l)}\in \omega_{l,\delta}\; |\; dist(\hat x^{(l)},J)<\eta_0
\delta\bigr\}}{\cal T}_\delta\bigl(\gamma_{ij}(V^{(l)}_\delta)\bigr)\longrightarrow 0
\qquad\hbox{strongly in}\quad L^2(\Omega_l).\leqno(3.21)$$ The part of ${\cal S}_\delta$ neighbour of the edge
$J$ and common to several plates is contained into the cylinder $\displaystyle {\bigl\{x\in \R^3\; |\; dist(x,J)<\eta_0
\delta\bigr\}}$. Thanks to the convergences $(3.21)$ its contribution in the limit problem is equal to zero. Then we can make 
$\delta$ tends to 0 in order to obtain $(3.19)$ with the displacement $V$. The limit of the right handside term of
$(3.10)$ is given by $(3.13)$. 
\vskip 1mm
\noindent  The set $H^1_{\Gamma_0}(\GS,\R^3)$ is dense in $H^1_{\rho,\Gamma_0}(\GS,\R^3)$ (Lemma B.3), which gives 
$(3.19)$ with any displacement of $D_E(\GS)$. \fin

\noindent{\bf 3.9. The inextensional displacement  $U_I$ or the problem of coupled bending plates }
\vskip 1mm
\noindent{\bf Theorem 3.10 : }{\it  The inextensional displacement $U_I$ is the solution of the variational problem   
$${E\over 3(1-\nu^2)}\displaystyle\sum_{l=1}^N\int_{\omega_l}\Bigl[(1-\nu){\partial^2 U^{(l)}_{I,3}\over
\partial x^{(l)}_\alpha\partial x^{(l)}_\beta}{\partial^2 V^{(l)}_3\over\partial x^{(l)}_\alpha\partial x^{(l)}_\beta} +\nu\Delta
U^{(l)}_{I,3}\Delta V^{(l)}_3 \Bigr]=\int_{\GS}f_I\cdot V,\qquad\forall V\in {\cal D}_I(\GS).\leqno(3.22)$$} 
\noindent The proof of Theorem 3.10 requires the next lemma.
\smallskip
\noindent{\bf Lemma 3.11 : }{\it  For any element $V\in {\cal D}_I(\GS)$, there exists a sequence of displacements
$\bigl(W_\delta\bigr)_{ 0<\delta\le \delta_0}$ such that
$$W_\delta\in H^1_{\Gamma_0}({\cal S}_{\delta},\R^3),\quad \hbox{and}\quad \left\{\eqalign{ {\cal
T}_\delta\bigl(W^{(l)}_\delta\bigr)&\longrightarrow  V^{(l)} \quad\hbox{strongly in}\quad L^2(\Omega_l,\R^3),\cr  
{\cal T}_\delta\bigl(\gamma_{\alpha\beta}(W^{(l)}_\delta)\bigr)&\longrightarrow-t^{(l)}_3{\partial^2
V^{(l)}_3\over\partial x^{(l)}_\alpha\partial x^{(l)}_\beta} \quad\hbox{strongly in}\quad L^2(\Omega_l),\cr 
{\cal T}_\delta\bigl(\gamma_{k3}(W^{(l)}_\delta)\bigr)&\longrightarrow 0\quad\hbox{strongly in}\quad L^2(\Omega_l).\cr}
\right.\leqno(3.23)$$}
\noindent{\bf Proof  : } See Annex E.\fin 
\noindent{\bf Proof of Theorem 3.10 : } Let $W$ be an element of ${\cal D}_I(\GS)$.  For any edge $J$ and any face
$\overline{\omega}_l$ containing $J$ we have
$$\hbox{\bf 1}_{\bigl\{\hat x^{(l)}\in \omega_{l,\delta}\; |\; dist(\hat x^{(l)},J)<\eta_0
\delta\bigr\}}{\cal T}_\delta\bigl(\gamma_{ij}(W_\delta^{(l)})\bigr)\longrightarrow 0
\qquad\hbox{strongly in}\quad L^2(\Omega_l).\leqno(3.24)$$  
\noindent  We take $W_\delta$ as  a test-displacement in $(3.10)$. We transform, by unfolding, the integrals on the
 plates  into integrals on the reference plates, then we divide by $2\delta$. We pass to the limit  (thanks to $(3.24)$ the
contribution  of the immediate junction neighborhoods tends to 0).  We obtain $(3.22)$ with the
test-displacement $V$.\fin
\noindent{\bf Remark 3.12 : } The problems $(3.19)$ and $(3.22)$ are coercive. It results that
the whole encountered sequences converges to  their limit.  We are going to show now that these convergences are strong.  We 
consider the formal displacement $U$ of the structure ${\cal S}_\delta$ defined in each plate by
$$U^{(l)}(x)=U^{(l)}_E(\widehat{x}^{(l)})+{1\over
\delta}U^{(l)}_I(\widehat{x}^{(l)})+\widehat{\nabla}U_I^{(l)}(\widehat{x}^{(l)})\land x^{(l)}_3\vec e^{(l)}_3,\qquad  x\in
\Omega_{l,\delta}.$$

\noindent Let $\hbox{\bf 1}_{{\cal J}^c_{\delta}}$  be the characteristic function of the complement  in ${\cal S}_\delta$ of the
union of the  edges neighborhoods $\displaystyle \bigcup_{J\in {\cal J}}\bigl\{x\in {\cal S}_\delta\; |\; dist(x, J)<\eta_0
\delta\bigr\}$. In the reference plate $\Omega_l$, we have the convergences 
$${\cal T}_\delta(\gamma_{ij}(u^{(l)}_\delta)){\cal T}_\delta(\hbox{\bf 1}_{{\cal J}^c_{\delta}})\rightharpoonup
0\quad\hbox{weakly in}\quad L^2(\Omega_l).$$ Hence
$$\left\{\eqalign{
&\sum_{l=1}^N\int_{\Omega_l}\sigma_{ij}(U^{(l)})\gamma_{ij}(U^{(l)})\le \liminf_{\delta\rightarrow 0}\sum_{l=1}^N
\int_{\Omega_l}{\cal T}_\delta(\sigma_{ij}(u^{(l)}_\delta)){\cal T}_\delta(\gamma_{ij}(u^{(l)}_\delta)){\cal T}_\delta(
\hbox{\bf 1}_{{\cal J}^c_{\delta}})\cr 
= &\liminf_{\delta\rightarrow 0}{1\over 2\delta}\int_{{\cal S}_\delta}\sigma_{ij}(u_\delta)\gamma_{ij}(u_\delta)\hbox{\bf
1}_{{\cal J}^c_{\delta}}\le \limsup_{\delta\rightarrow 0}{1\over 2\delta}\int_{{\cal
S}_\delta}\sigma_{ij}(u_\delta)\gamma_{ij} (u_\delta)\cr 
=&\limsup_{\delta\rightarrow 0}{1\over2\delta}\int_{{\cal S}_\delta}F_\delta\cdot u_\delta =\lim_{\delta\rightarrow
0}{1\over2\delta}\int_{{\cal S}_\delta}F_\delta\cdot u_\delta=\int_\GS f_E\cdot U_E+\int_\GS f_I\cdot
U_I\cr}\right.$$
\noindent The first term of $(3.24)$ is the sum of the left handside members of $(3.19)$ and $(3.22)$. Hence the
above inequalities are equalities. Besides
$$\eqalign{
\liminf_{\delta\rightarrow 0}{1\over 2\delta}\int_{{\cal S}_\delta}\sigma_{ij}(u_\delta)\gamma_{ij}(u_\delta) 1_{{\cal
J}^c_{\delta}} &\le\liminf_{\delta\rightarrow 0}{1\over 2\delta}\int_{{\cal
S}_\delta}\sigma_{ij}(u_\delta)\gamma_{ij}(u_\delta)\le \limsup_{\delta\rightarrow 0}{1\over 2\delta}\int_{{\cal
S}_\delta}\sigma_{ij}(u_\delta)\gamma_{ij} (u_\delta)\cr
\hbox{and}\quad \liminf_{\delta\rightarrow 0}{1\over 2\delta}\int_{{\cal S}_\delta}\sigma_{ij}(u_\delta)\gamma_{ij}(
u_\delta)1_{{\cal J}^c_{\delta}} &
\le\limsup_{\delta\rightarrow 0}{1\over 2\delta}\int_{{\cal S}_\delta}
\sigma_{ij}(u_\delta)\gamma_{ij}(u_\delta)1_{{\cal J}^c_{\delta}}\le \limsup_{\delta\rightarrow 0}{1\over
2\delta}\int_{{\cal S}_\delta}\sigma_{ij}(u_\delta)\gamma_{ij} (u_\delta)}$$
\noindent We deduce that $\displaystyle \lim_{\delta\rightarrow 0}{1\over 2\delta}\int_{{\cal
S}_\delta}\sigma_{ij}(u_\delta)\gamma_{ij}(u_\delta)\bigl(1- \hbox{\bf 1}_{{\cal J}^c_{\delta}}\bigr)=0$. All the sequences of
the unfolded of the strained tensor components strongly converge in $L^2(\Omega_l)$. We have also the strong  convergences
$$\left\{\eqalign{
\delta {\cal U}_\delta&\longrightarrow U_I,\quad\delta{\cal R}_\delta \longrightarrow  {\cal R}\quad\hbox{strongly in}\quad
H^1(\GS,\R^3),\cr 
\delta U_{I,\delta}&\longrightarrow U_I\qquad\hbox{strongly in}\quad D_I(\GS),\cr 
U_{E,\delta}&\longrightarrow  U_E\qquad\hbox{strongly in}\quad D_E(\GS).\cr }\right.$$
\noindent{\bf 3.10. Complements }
\vskip 1mm
\noindent The orthogonal condition  $(3.12)$ requires an explanation . First, for any function
$\phi\in H^1(\omega_l)$ equal to zero on the edges, the displacement $\Phi$ defined by
$$\Phi^{(l)}=\phi\vec e_3^{(l)},\quad \hbox{in }\omega_l,\; \hbox{and by } 0 \; \hbox{in the other faces,}$$ belongs to
 $D_I(\GS)$.  We deduce that the function  $f_{E,3}^{(l)}$ is orthogonal to $\phi$ and then, by density of these test-functions
in $L^2(\omega_l)$, we get
$$\forall l\in\{1,\ldots,N\},\qquad f_{E,3}^{(l)}=0.\leqno(3.25)$$
\noindent Let  $D_{I,0}(\GS)$ be the space of the inextensional displacements equal to zero on the edges belonging to
${\cal J}$ and let $\bigl(D_{I,0}(\GS)\bigr)^\perp$ be its orthogonal in $D_I(\GS)$ for the inner product $<\cdot,\cdot>_\rho$.
The subset $\bigl(D_{I,0}(\GS)\bigr)^\perp$ is of finite dimension. The condition $(3.12)$ is then equivalent to
$$\forall V\in \bigl(D_{I,0}(\GS)\bigr)^\perp,\qquad \int_\GS f_E\cdot V=0.\leqno(3.26)$$  This last condition results in a finite
number of equalities related to  the means in the faces $\overline{\omega}_l$ of the components $f^{(l)}_{E,\alpha}$ of $f_E$.
\bigskip
\noindent{\bf 4. Annexes }
\vskip 1mm
\noindent{\bf  4.1 Annex A. Proof of Theorem 2.3}
\vskip 1mm
\noindent  The proof of Theorem 2.3  is based  on Lemma 2.3 in [4] and on Lemma 4.1. 

\noindent We denote
$$\omega_\eta=\bigl\{\widehat{x}\in \R^2\; |\; dist(\widehat{x},\omega)<\eta\bigr\},\qquad \eta>0.$$
\noindent{\bf Lemma 4.1  : }{\it There exist  $R>0$ and $\delta_0>0$, depending only on $\omega$, such that for any
$\delta\in]0,\delta_0]$, $\omega_{2\delta}$ is covered by a family of open sets, of diameter less than
 $R\delta$, star-shaped with respect to a disc of radius ${\delta/ 2}$ and such that any point of $\omega_{2\delta}$
belongs to a finite number (independent of $\delta$) of open sets of that family.} 
\vskip 1mm
\noindent{\bf Proof : } The open set ${\cal A}_{pq}=](p-1/2)\delta,(p+3/2)\delta[\times](q-1/2)\delta,(q+3/2)\delta[$,
$(p,q)\in \Z^2$, has a diameter of $2\sqrt 2\delta$ and is star-shaped with respect to the disc of center
$\bigl((p+1/2)\delta,(q+1/2)\delta)\bigr)$ and of radius ${\delta/ 2}$. Let ${\cal I}_{\delta}$ be the set of the pairs
$(p,q)$ of $\Z^2$ such that ${\cal A}_{pq}\subset\omega$. The distance between the boundary of $\omega$ and 
$\bigcup_{(p,q)\in {\cal I}_{\delta}}{\cal A}_{pq}$ is less than $3\delta$.
\medskip
\noindent Let us proceed now to the covering of the neighborhood of the boundary of $\omega$. 

\noindent The boundary of $\omega$ is lipschitzian. Hence there exist constants $A$,  $B$, $C$, $M$ 
strictly positive, a finite number $N$ of local coordinate systems $(x_{1r},x_{2r})$ in $(O_r;\vec e_{1r},\vec
e_{2r})$ and maps $f_r\,: \, [-A,A]\longrightarrow \R$, Lipschitz continuous with ratio $M$, $1\le r\le N$, such that  
$$\left\{\eqalign{
& \partial \omega=\displaystyle\bigcup_{r=1}^N\Bigl\{(x_{1r},x_{2r})\; |\;  x_{2r}=f_r(x_{1r}),\quad x_{1r}\in]-A,A[\Bigr\},\cr
&\bigl\{x\in\omega\, |\, dist(x,\partial\omega)<C\bigr\}\i \bigcup_{r=1}^N\Bigl\{(x_{1r},x_{2r})\ \;
|\; f_r(x_{1r})<x_{2r}<f_r(x_{1r})+B,\; |x_{1r}|\le A\Bigr\}\subset\omega,\cr
&\omega_C\setminus\omega\i \bigcup_{r=1}^N\Bigl\{(x_{1r},x_{2r})\; |\; 
f_r(x_{1r})-B<x_{2r}<f_r(x_{1r}),\;|x_{1r}|\le A\Bigr\}\subset\R^2\setminus\omega.\cr}\right.$$ Through the use of easy
geometrical arguments we show that if  $4\delta\le
\inf\{C,{B/ \sqrt{1+M^2}}\}$, we have
$$\eqalign{ &\omega_{2\delta}\setminus \omega\i\bigcup_{r=1}^N\Bigl\{(x_{1r},x_{2r}) \;|\;
f_r(x_{1r})-2\delta\sqrt{1+M^2} <x_{2r}<f_r(x_{1r}),\quad |x_{1r}|\le A\Bigr\},\cr &\bigl\{x\in\omega\, |\,
dist(x,\partial\omega)<4\delta\bigr\} 
\i  \bigcup_{r=1}^N\Bigl\{(x_{1r},x_{2r}) \;|\; f_r(x_{1r})<x_{2r}<f_r(x_{1r})+4\delta\sqrt{1+M^2} ,\,
|x_{1r}|\le A\Bigr\}.\cr}$$
\noindent For any $a\in]-A,A-2\delta[$, the domains   
$${\cal B}_{\delta, a, r}=\bigl\{(x_{1r},x_{2r})\; |\; f_r(x_{1r})-(6M+2)\delta<x_{2r}<f_r(x_{1r})+(6M+2)\delta,\enskip
x_{1r}\in]a,a+2\delta[\bigr\}$$ and
${\cal B}_{\delta, a, r}\cap \omega$ are star-shaped with respect to the disc of center
$(a+\delta,f_r(a)+(3M+1)\delta)$ and of radius ${\delta / 2}$. These open sets have a diameter less than 
$6(3M+1)\delta=R\delta$. 

\noindent For $0<\delta\le \delta_0=\inf\bigl\{{B/(6M+2)}, A/2, C/4\bigr\}$ the open sets  ${\cal A}_{pq}$
$\bigl((p,q)\in {\cal I}_{\delta}\bigr)$, ${\cal B}_{\delta, a_p, r}$, where $a_p={p\delta}\in[-A,A[$ ($p\in\Z$), ${\cal
B}_{\delta, -A, r}$ and ${\cal B}_{\delta, A-2\delta, r}$ ($r\in \{1,\ldots, N\}$) cover 
$\omega_{2\delta}$; their diameter is less than $R\delta$  and they are star-shaped with respect to a disc of radius
$\delta/2$. Any point of $\omega_{2\delta}$ belongs to a finite number (depending only on $\omega$) of open sets of
that family.\fin

\noindent We denote $\bigl\{\omega^{'}_{\delta,n}\bigr\}_{n\in N_\delta}$ the covering of $\omega_{2\delta}$  obtained in
Lemma A.1  and $\bigl\{\omega_{\delta,n}\bigr\}_{n\in N_\delta}$ the covering of 
$\omega$ defined by $\omega_{\delta,n}=\omega^{'}_{\delta,n}\cap \omega$, $n\in N_\delta$.
$$\qquad\vbox{\centerline{\dessin 195.7mm by 257.8mm (Fig2 scaled 300)}} $$
\centerline{Figure 2. The domain ${\cal B}_{\delta, a, r}$}
\noindent{\bf Proof of Theorem 2.3 : } The open set $\omega_{\delta,n}$ is star-shaped with respect to a disc of center
$A_n$ and of radius ${\delta/2}$. We put ${\cal O}_{\delta,n}=\omega_{\delta,n}\times ]-\delta,\delta[\i \Omega_\delta$,
$n\in N_\delta$. The domain ${\cal O}_{\delta,n}$ has a diameter less than $(R+2)\delta$, and is star-shaped with respect to a
ball of center $A_n$ and of radius ${\delta/ 2}$. From Lemma 2.3 of [4], there exists a rigid displacement $r_n$ such that
$${\cal D}(u-r_n,{\cal O}_{\delta,n})+{1\over \delta^2}||u-r_n||^2_{L^2({\cal O}_{\delta, n} , \R^3)}\le C{\cal E}(u,{\cal
O}_{\delta,n}),\quad r_n(x)=a_n+b_n\land \overrightarrow{A_nx}, \quad (a_n,b_n)\in\R^3.\leqno (4.1)$$ 
\noindent The constant  depends only on $R$. 

\noindent We calculate the mean of $(u-r_n)(x)$ and of $x_3\vec e_3\land (u-r_n)(x)$ on the intervals
$\{\widehat{ x}\}\times]-\delta,\delta[$, $\widehat{x}\in \omega_{\delta,n}$, then we integrate on $\omega_{n,\delta}$ the
inequalities we have obtained. Thanks to $ (4.1)$, we have
$$\int_{\omega_{\delta,n}}|{\cal U}(\widehat{ x})-a_n-b_n\land\overrightarrow{A_n\widehat{x}}|^2d\widehat{ x} \le C\delta
{\cal E}(u,{\cal O}_{\delta, n}),\qquad \int_{\omega_{\delta,n}}|{\cal R}_\alpha(\widehat {x})-b_{\alpha,n}|^2d\widehat{ x}\le
{C\over\delta}{\cal E}(u,{\cal O}_{\delta,n}).\leqno (4.2)$$  In  $ (4.1)$ we eliminate the displacement $r_n$ thanks
to the estimations $(4.2)$. Hence we have  $||u-U_e||^2_{L^2({\cal O}_{\delta,n} , \R^3)}\le C\delta^2{\cal E}(u,{\cal
O}_{\delta,n})$, then we add these inequalities and we obtain  
$$ ||u-U_e||^2_{L^2(\Omega_\delta,\R^3)}\le C\delta^2{\cal E}(u,\Omega_\delta).$$
\noindent Both components of e.p.d. $U_e$ belong to $H^1(\omega,\R^3)$. We calculate the mean of the gradient of
$(u-r_n)(x)$, then the mean of $ x_3\vec e_3\land \nabla(u-r_n)(x)$ on the intervals $\{\widehat{
x}\}\times]-\delta,\delta[$, $\widehat{x}\in
\omega_{\delta,n}$. Using $(4.1)$ we obtain 
$$\Bigr\|{\partial{\cal U}\over\partial x_\alpha}-b_n\land\vec e_\alpha\Bigr\|^2_{L^2(\omega_{\delta,n},\R^3)}
+\delta^2\Bigr\|{\partial{\cal R}\over\partial x_\alpha}\Bigr\|^2_{L^2(\omega_{\delta,n},\R^2)}\le {C\over \delta}{\cal
E}(u,{\cal O}_{\delta,n}),\leqno (4.3)$$ hence, after elimination of $b_n$ in the first inequality,
$$\left\{\eqalign{ &\Bigr\|{\partial{\cal U}_1\over\partial x_1}\Bigr\|^2_{L^2(\omega_{\delta,n})}+\Bigr\|{\partial{\cal
U}_2\over\partial x_2}\Bigr\|^2_{L^2(\omega_{\delta,n})}+\Bigr\|{\partial{\cal U}_1\over\partial x_2}+{\partial{\cal
U}_2\over\partial x_1}\Bigr\|^2_{L^2(\omega_{\delta,n}) }\cr +&\Bigr\|{\partial{\cal U}_3\over\partial x_1}+{\cal
R}_2\Bigr\|^2_{L^2(\omega_{\delta,n})}+\Bigr\|{\partial{\cal U}_3
\over\partial x_2}-{\cal R}_1\Bigr\|^2_{L^2(\omega_{\delta,n})}\le {C\over \delta}{\cal E}(u,{\cal
O}_{\delta,n})\cr}\right.\leqno (4.4)$$ From $ (4.3)$ and $ (4.4)$ we deduce the estimate of ${\cal E}(U_e,{\cal
O}_{\delta,n})$
$${\cal E}(U_e,{\cal O}_{\delta,n})\le C{\cal E}(u,{\cal O}_{\delta,n})\qquad \hbox{hence}\qquad {\cal
E}(U_e,\Omega_\delta) \le C{\cal E}(u,\Omega_\delta).$$ From $ (4.3)$, $ (4.4)$, $(4.1)$ and
after elimination of the gradient of $r_n$ we also deduce
$${\cal D}(u-U_e,{\cal O}_{\delta,n})\le C{\cal E}(u,{\cal O}_{\delta,n}),\quad\hbox{hence}\quad {\cal
D}(u-U_e,\Omega_\delta)\le C{\cal E}(u,\Omega_\delta).$$  Theorem 2.3  is proved.\fin
\vskip 1mm
\noindent{\bf  4.2 Annex B. About the second decomposition of a plate displacement}
\vskip 1mm
\noindent{\bf  4.2.a Extension of a plate displacement}
\vskip 1mm
\noindent Let $\omega$ be a polygonal bounded domain in $\R^2$.  The boundary of $\omega$ is made of a finite number of
segments. Let ${\cal C}$ be a connected component of $\partial \omega$. There exists $\delta^{'}_0>0$ such that
for any  $\delta\in]0,\delta^{'}_0]$ the domains
$${\cal C}_\delta=\Bigl\{x\in \R^3\; |\; dist(x,{\cal C})<\delta\Bigr\},\quad\hbox{and}\quad  {\cal C}_{3\delta}=\Bigl\{x\in
\R^3\; |\; dist(x,{\cal C})< 3\delta\Bigr\}$$ are  rods structures. Then there exists $\mu^{'}_0>0$ such that all the balls
  centered in a vertex of  ${\cal C}$, and of radius $3\eta_0^{'}\delta$ contain the junctions of the rods belonging to
${\cal C}_{3\delta}$.

\noindent We recall that for any $\delta\in]0,\delta{'}_0]$, there exists an extension  operator, linear
and  continuous, $P_\delta$ from $H^1({\cal C}_\delta)$ into $H^1({\cal C}_{3\delta})$ such that for any $\phi\in H^1({\cal
C}_\delta)$,
$$P_\delta(\phi)_{|{\cal C}_\delta}=\phi,\qquad  ||P_\delta(\phi)||_{L^2({\cal C}_{3\delta})}+\delta||\nabla P_\delta(\phi)
||_{[L^2({\cal C}_{3\delta})]^3} \le  C\Bigl\{||\phi||_{L^2({\cal C}_{\delta})}+\delta||\nabla\phi||_{[L^2({\cal C}_\delta)]^3}
\Bigr\}.\leqno(4.5)$$ The constant does not depend on $\delta$.

\noindent{\bf Proof of Lemma 2.8 : }  We begin with extending $u$ in the neighborhood of a connected component of
$\partial \omega$. 

\noindent Let ${\cal C}$ be a  connected component of $\partial \omega$. The restriction of $u$ to ${\cal
C}_\delta$ is a displacement belonging to $H^1({\cal C}_\delta,\R^3)$. Hence there exists an elementary displacement
of a rods structure (e.d.r.s.) $U_{e,R}$ (see [5])  which coincides with a rigid displacement in each set
$B(A,3\eta^{'}_0\delta)\cap {\cal C}_\delta$ where $A$ is a vertex of ${\cal C}$ and which verifies
$${\cal E}(U_{e,R},{\cal C}_\delta)+{\cal D}(u-U_{e,R},{\cal C}_\delta)+{1\over \delta^2}||u-U_{e,R}||^2_{L^2( {\cal
C}_\delta,\R^3)}\le C{\cal E}(u,{\cal C}_\delta).\leqno(4.6)$$ The displacement  $U_{e,R}$ is also an  e.d.r.s.  of
${\cal C}_{3\delta}$ and ${\cal E}(U_{e,R},{\cal C}_{3\delta})\le C{\cal E}(u,{\cal C}_\delta)$. The displacement
$$u^{'}_{{\cal C}_{3\delta}}=U_{e,R}+P_\delta(u-U_{e,R})$$ is an extension of $u$ to the set ${\cal C}_{3\delta}$.
From $(4.6)$ we have the following inequalities:  
$${\cal E}(u^{'}_{{\cal C}_{3\delta}},{\cal C}_{3\delta})\le 2\bigl\{{\cal E}(U_{e,R},{\cal C}_{3\delta})+{\cal
E}(P_\delta(u-U_{e,R}),{\cal C}_{3\delta})\bigr\}\le  C{\cal E}(u,{\cal C}_\delta)\leqno(4.7)$$  
In the same way we build an extension of $u$ in the neighborhood of the other connected components of
$\partial\omega$. The extension ${\cal P}_\delta(u)$  is then the displacement which coincides with $u$ in
$\Omega^{'}_\delta$ and which is equal to one of the previous extensions in
$\widetilde{\Omega}_\delta\setminus
\Omega^{'}_\delta$. The  estimates $(2.14)$ are the immediate consequences of the inequalities 
$(4.7)$ obtained in the neighborhood of each connected components of 
$\omega$.\fin
\noindent{\bf  Remark 4.2 : } If one of the edges of the boundary of $\omega$ is
fixed we can take an e.d.r.s. with its two components  equal to zero on this edge without modifying the
estimates $(4.6)$ and then extend $u$ beyond this edge by $0$.\fin
\noindent{\bf Remark 4.3 : } We also can construct an extension operator $P_\delta$  when  $\omega$  is of  lipschitzian 
boundary with the help of a few changes.\fin

\noindent{\bf 4.2.b  Modification of an e.p.d. in the neighborhood of an edge.}

\noindent Let $J$ be an edge contained in the face $\overline{\omega}$, $J_\delta$ the rod
$$J_\delta=\bigl\{x\in \R^3\; |\; dist(x,J)<\delta\bigr\}\i \Omega{'}_\delta,$$ and $u$  a displacement of the plate
$\widetilde{\Omega}_\delta$. Without being detrimental to the general case we can suppose that the edge's direction is
$\vec e_1$ and that one of these extremities is the chosen origin on the face, so that $J$ is identified with the segment
$[0,L]\times \{0\}$ where $L$ is the edge's length.

\noindent The restriction of the displacement $u$ to the rod $J_\delta$ can be decomposed into the sum of the
elementary rod displacement (e.r.d.) $U_{e,R}$ of components ${\cal U}_R$ and ${\cal R}_R$ and of a residual displacement. We
choose an e.r.d. $U_{e,R}$  coinciding with a  rigid displacement in the balls centered in the extremities of $J$ and of radius
$\eta_0\delta$ (see [5]). We have
$$U_{e,R}(x)={\cal U}_{R}(x_1)+{\cal R}_{R}(x_1)\land\bigl(x_2\vec e_2 +x_3\vec e_3\bigr).$$ We know (see [4] and
[5]) that the components ${\cal U}_R$ and ${\cal R}_R$ of $U_{e,R}$ belong to
$H^1(J,\R^3)$, and verify
$$\left\{\eqalign{ &
\delta^2\Bigl\|{d{\cal R}_{R}\over dx_1}\Bigr\|^2_{L^2(]0,L[,\R^3)} + \Bigl\|{d{\cal U}_{R}\over dx_1} -{\cal R}_{R}\land\vec
e_1\Bigr\|^2_{L^2(]0,L[,\R^3)}\le {C\over \delta^2}{\cal E}(u, J_\delta)\cr 
&{\cal E}(U_{e,R},J_\delta)+{\cal D}(u-U_{e,R},J_\delta)+{1\over
\delta^2}||u-U_{e,R}||^2_{L^2(J_\delta,\R^3)}\le C{\cal E}(u, J_\delta)\cr}\right.\leqno(4.8)$$ The functions
${\cal U}_R$ and ${\cal R}_R$  are extended into functions belonging to $H^1_{loc}(\R,\R^3)$ (by construction ${\cal
R}_R$ is constant and ${\cal U}_R$ is linear in a neighborhood of the extremities of  $J$). These extensions are then identified
with elements belonging to $H^1_{loc}(\R^2,\R)$ depending only on the variable $x_1$.

\noindent Let 
$$J^{'}_{\delta}=\Bigl\{x\in \Omega^{'}_\delta\; |\; dist(\widehat{x}, J)< (2\eta_0+1) \delta\;\Bigr\}$$ be the rod and
$\widehat{J}_\delta$ the neighborhood of $J$ in $\omega_\delta$,
$$\widehat{J}_{\delta}=\Bigl\{x\in \omega_\delta\; |\; dist(\widehat{x}, J)< 2\eta_0\delta\;\Bigr\}$$ From
the estimates
$(4.8)$ of the restriction of $u-U_{e,R}$ to $J_\delta$ we deduce the following estimates of the restriction of
$u-U_{e,R}$
to $J^{'}_\delta$: 
$${\cal D}(u-U_{e,R},J^{'}_\delta)+{1\over \delta^2}||u-U_{e,R}||^2_{L^2(J^{'}_\delta,\R^3)}\le C{\cal E}(u,
J^{'}_\delta).\leqno(4.9)$$ The constant depends on $L$ and   $\eta_0$.

\noindent The displacement $u$ of the plate $\Omega^{'}_\delta$  is decomposed into the sum of an elementary plate
displacement $U_{e,P}$, given by $(2.15)$, and of a residual  displacement,
$$U_{e, P}(x)={\cal U}_P(\widehat{x})+{\cal R}_P(\widehat{x})\land x_3\vec e_3,\qquad x\in \Omega^{'}_\delta.$$
Besides the inequalities $(2.16)$, we also have
$$\left\{\eqalign{
\delta^3||\nabla{\cal R}_P||^2_{L^2(\widehat{J}_{\delta},\R^6)}+\delta\Bigl\|{\partial {\cal U}_P\over \partial
x_\alpha}- {\cal R}_P\land \vec e_\alpha\Bigr\|^2_{L^2(\widehat{J}_{\delta},\R^3)}\le
C{\cal E}(u, J^{'}_\delta),\cr
{\cal D}(u-U_{e,P},J^{'}_\delta)+{1\over \delta^2}||u-U_{e,P}||^2_{L^2(J^{'}_\delta,\R^3)}\le C{\cal E}(u,
J^{'}_\delta),\cr}\right.\leqno(4.10)$$ This allows us to compare the different elementary displacements. 
We obtain
$$||U_{e,R}-{\cal U}_P||^2_{L^2(\widehat{J}_{\delta},\R^3)}+\delta^2||{\cal R}_{R}-{\cal R}_{P}||^2_{
L^2(\widehat{J}_{\delta},\R^3)}\le C \delta{\cal E}(u,J^{'}_\delta).\leqno(4.11)$$ The estimate of $||{\cal R}_{R}-{\cal
R}_{P}||^2_{ L^2(\widehat{J}_{\delta},\R^3)}$ follows from the nullity of ${\cal R}_{R}-{\cal R}_{P}$ on $J$.
\vskip 1mm
\noindent We are now going to modify the  e.r.d. $U_{e, P}$ in the neighborhood of $J$.

\noindent We consider a function  $m$ belonging to ${\cal C}^\infty(\R^+, [0,1])$ such that  
$$ m(t)=1\qquad\forall t\ge 2,\qquad\qquad m(t) =0\qquad\forall t\le 1,\qquad\qquad |m^{'}(t)| \le 2\qquad\forall
t\in\R.\leqno(4.12)$$ We define the components, ${\cal U}^{'}$ and ${\cal R}^{'}$, of a new e.r.d. $U^{'}_e$ by
$$\left\{\eqalign{ 
{\cal U}^{'}(\widehat{x})&=U_{e,R}(\widehat{x})\Bigl(1-m\Bigl({dist(\widehat{ x}, J)\over\eta_0
\delta}\Bigr)\Bigr)+{\cal U}_P(\widehat{x})m\Bigl({dist(\widehat{ x}, J)\over \eta_0\delta}\Bigr),\cr 
{\cal R}^{'}(\widehat{x})&={\cal R}_{R }(\widehat{x})\Bigl(1-m\Bigl({dist(\widehat{ x}, J)\over
\eta_0\delta}\Bigr)\Bigr)+{\cal R}_{P}(\widehat{x})m\Bigl({dist(\widehat{ x}, J)\over\eta_0\delta}\Bigr),\qquad
\widehat{x}\in\omega_{\delta},\cr
U^{'}_{e }(x)&={\cal U}^{'}(\widehat{x})+{\cal R}^{'} (\widehat{x})\land x_3\vec e_3,\qquad x\in \Omega{'}_\delta.}
\right.\leqno(4.13)$$ 
\noindent Hence we have by construction of $U^{'}_e$,
$$\eqalign{
&\hbox{if}\enskip x\in\Omega^{'}_\delta \enskip\hbox{and} \enskip dist(\widehat{ x}, J)
<\eta_0\delta\quad\hbox{then}\quad \; U^{'}_e(x)= U_{e,R}(x)\cr
&\hbox{if}\enskip x\in\Omega^{'}_\delta \enskip\hbox{and} \enskip dist(\widehat{ x}, J)
>2\eta_0\delta\quad\hbox{then}\quad U^{'}_e(x)= U_{e,P}(x).\cr}$$
\noindent Thanks to $(4.8)$, $(4.10)$,  $(4.12)$ and $(2.16)$ the e.p.d. $U^{'}_e$  verifies
$$\eqalign{ 
&\delta^3||\nabla{\cal R}^{'}||^2_{L^2(\omega_\delta,\R^6)}+\Bigl\|{\partial {\cal U}^{'}\over \partial x_\alpha}
-{\cal R}^{'}\land \vec e_\alpha\Bigr\|^2_{L^2(\omega_\delta,\R^3)}+{\cal E}(U^{'}_e,\Omega{'}_\delta)\le C{\cal
E}(u,\Omega{'}_\delta),\cr  &{\cal D}(u-U^{'}_e, \Omega^{'}_\delta)\le C{\cal E}(u, \Omega^{'}_\delta),\qquad
||u-U^{'}_e||^2_{L^2( \Omega{'}_\delta, \R^3)}\le C\delta^2{\cal E}(u, \Omega^{'}_\delta).\cr}$$ The
constants depend only on $\omega$, $J$ and $\eta_0$. 
\vskip 1mm
\noindent{\bf 4.3 Annex C. About the spaces $H^1_{\rho,\Gamma_0}(\GS,\R^3)$ and $D_E(\GS)$ }
\vskip 1mm
\noindent For any $\theta_0\in]0,\pi[$ and any $r>0$, we denote $\displaystyle C_{r,\theta_0}=\Bigl\{(x,y)\in B\bigl(O;r\bigr)\;
|\;  0<\theta<\theta_0\Bigr\}$,  $J_0$ the  segment  of origin $O$ and of extremity $A=(0,1)$ and we denote $J_{\theta_0}$
the  segment of origin $O$ and extremity $B=(\cos(\theta_0),\,\sin(\theta_0))$. 
\vskip 1mm
\noindent We denote $r=dist(x, O)$, $x\in \R^2$. 
\vskip 1mm
\noindent{\bf Lemma 4.3 : }{\it Let $\phi$ belong to $ H^1(C_{1,\theta_0})$, for any $\alpha\in]0,1]$, we have 
$$\int_{C_{1,\theta_0}}|\phi|^2r^{\alpha-2}\le {4\over
\alpha}||\phi||^2_{L^2(C_{1,\theta_0})}+{2\over\alpha^2}||\nabla\phi||^2_{[L^2(C_{1,\theta_0})]^2}\leqno(4.14)$$}
\noindent{\bf Proof : } We recall that  for any $u\in H^1(0,L)$ and  for any $\alpha\in]0,1]$, we have
$$\int_0^L|u(t)|^2  t^{\alpha-1}dt\le {4\over \alpha L^{2-\alpha}}\int_0^L |u(t)|^2tdt+{2\over
\alpha^2L^\alpha}\int_0^L|u{'}(t)|^2tdt\leqno(4.15)$$
Let us take $\phi\in{\cal C}^\infty(\overline{C}_{1,\theta_0})$. We apply the inequality 4.4 to the restriction of $\phi$
to a radius coming from the origin and contained in  $C_{1,\theta_0}$. This gives
$$\eqalign{
\int_0^1{|\phi\bigl(r\cos(\theta),r\sin(\theta)\bigr)|^2\over r^{1-\alpha}}dr &\le
{4\over\alpha}\int_0^1|\phi\bigl(r\cos(\theta), r\sin(\theta)\bigr)|^2rdr
+{2\over\alpha^2}\int_0^1\Bigl|{\partial\phi\over\partial r}\bigl(r\cos(\theta),r\sin(\theta)\bigr)\Bigr|^2rdr.\cr}$$ We then
integrate with respect to  $\theta$ between $0$ and $\theta_0$ and we obtain $(4.14)$. The density of ${\cal
C}^\infty(\overline{C}_{1,\theta_0})$ into $H^1(C_{1,\theta_0})$ gives the inequality for any function of the space
$H^1(C_{1,\theta_0})$.\fin  
\noindent{\bf Lemma 4.4 : }{\it  Let $u$ be in  $H^{1/ 2}(J_0)$, $v$ be in  $H^{1/ 2}(J_{\theta_0})$ and $\alpha\in]0,2]$.
There exists a function $w$ belonging to $H^1_{loc}(C_{1,\theta_0})\cap L^2(C_{1,\theta_0})$ such that
$$\eqalign{ 
& w_{|J_0}=u,\qquad w_{|J_{\theta_0}}=v,\qquad \hbox{and}\cr
&||w||^2_{L^2(C_{1,\theta_0})}+\int_{C_{1,\theta_0}}|\nabla w|^2r^{\alpha}\le{C\over\alpha^2}
\bigl\{||u||^2_{H^{1/ 2}(J_0)} +||v||^2_{H^{1/2}(J_{\theta_0})}\bigr\}.\cr}$$ The constant depends only on 
$\theta_0$. Moreover $w$ belongs to $W^{1,p}(C_{1,\theta_0})$ for any $p$, such that $\displaystyle 1\le p<{4\over
2+\alpha}$.}   
\vskip 1mm  
\noindent {\bf Proof : } We denote $J^{'}_0$ (resp. $J^{'}_{\theta_0}$), the   segment of same direction as
$J_0$ (resp. $J_{\theta_0}$) and of length $\displaystyle \tan\bigl({\theta_0 / 2}\bigr)$. The function $u$ (resp. $v$)
extends by reflexion into an element still denoted $u$ (resp. $v$) belonging to $H^{1/2}(J^{'}_0)$ (resp.
$H^{1/2}(J^{'}_{\theta_0})$).

\noindent Let $\tilde u$ and $\tilde v$ be the functions belonging to $H^{1/2}(J^{'}_0\cup J^{'}_{\theta_0})$
defined by 
$$\eqalign{ 
&\tilde u_{|J{'}_0}=u,\qquad \tilde u_{|J^{'}_{\theta_0}}(t\vec e_{\theta_0})=u(t\vec e_0),\cr &\tilde
v_{|J^{'}_0}(t\vec e_0)=v(t\vec e_{\theta_0}),\qquad \tilde v_{|J^{'}_{\theta_0}}=v,\cr}\qquad t\in ]0,\tan ({\theta_0 /
2})[$$ where $\vec e_\theta=\cos(\theta)\vec e_1+\sin(\theta)\vec e_2,\; \theta\in [0,\theta_0]$. There exists a continuous
lifting operator from $H^{1/2}(J^{'}_0\cup J^{'}_{\theta_0})$ into  $H^1(C_{\tan(\theta_0/2),\theta_0})$. Let  $U$ (resp.
$V$) be the lifting of $\tilde u$ (resp. $\tilde v$).  In the triangle $T_{\theta_0}$ of vertexes $\bigl(0,0\bigr)$,
$\bigl(0,\tan(\theta_0/2)\bigr)$ and  $\tan(\theta_0/2)\bigl(\cos(\theta_0), \sin(\theta_0)\bigr)$,
containing $C_{1,\theta_0}$ and contained in $C_{\tan(\theta_0/2),\theta_0}$, we define $w$ by
$$w(x_1,x_2)= U(x_1,x_2){x_1\sin(\theta_0)-x_2\cos(\theta_0)\over x_2(1-\cos(\theta_0))+
x_1\sin(\theta_0)}+V(x_1,x_2){x_2\over x_2(1-\cos(\theta_0))+x_1\sin(\theta_0)}.$$ In the above expression the
coefficients of $U(x_1,x_2)$ and $V(x_1,x_2)$ are  barycentric coordinates of  point $(x_1,x_2)$ belonging to $T_{\theta_0}$.
By construction  we have $w_{|J_0}=u$  and $w_{|J_{\theta_0}}=v$. The function $w$   belongs to $ L^2(T_{\theta_0})$ and
$$||w||_{L^2(T_{\theta_0})}\le ||U||_{L^2(C_{\tan(\theta_0/2),\theta_0})}+||V||_{L^2(C_{\tan(\theta_0/2),\theta_0})}\le
C\bigl\{||u||^2_{H^{1/ 2}(J_0)} +||v||^2_{H^{1/2}(J_{\theta_0})}\bigr\}.$$ We then calculate the  partial derivatives of $w$ and
we conclude that  $w$ belongs to $H^1_{loc}(T_{\theta_0})$. Moreover we have 
$$|\nabla w(x_1,x_2)|\le C\bigl\{|\nabla U(x_1,x_2)|+|\nabla V(x_1,x_2)|+|(U-V)(x_1,x_2)| r^{-1}\bigr\}\qquad (x_1,x_2)\in
C_{\tan(\theta_0/2),\theta_0}$$ The constant depends on $\theta_0$.  Thanks to the inequality of Lemma 4.3, we have
$$\eqalign{
\int_{C_{1,\theta_0}}|U-V|^2 r^{\alpha-2} &\le {4\over \alpha}||U-V||^2_{L^2(C_{1,\theta_0})}
+{2\over\alpha^2}||\nabla(U-V)||^2_{[L^2(C_{1,\theta_0})]^2}\cr
\Longrightarrow\enskip\int_{C_{1,\theta_0}}|\nabla w|^2r^\alpha  &\le{C\over
\alpha^2}\bigl\{||U||^2_{H^1(C_{\tan(\theta_0/2),\theta_0})} +||V||^2_{H^1(C_{\tan(\theta_0/2),\theta_0})}\bigr\}.\cr}$$
Eventually we obtain the estimate of Lemma 4.4. Moreover we have
$$\eqalign{
\int_{C_{1,\theta_0}}|\nabla w|^p\le \Bigl\{\int_{C_{1,\theta_0}}|\nabla w|^2r^\alpha\Bigr\}^{p\over
2}\Bigr\{\int_{C_{1,\theta_0}}r^{-{\alpha p\over 2-p}}\Bigr\}^{2-p\over 2}\cr}$$ Hence $w$ belongs to 
$W^{1,p}(C_{1,\theta_0})$ if $\displaystyle {\alpha p\over 2-p}<2$. That is to say for $\displaystyle 1\le p<{4\over
2+\alpha}$.\fin
\noindent{\bf Corollary : } If $\alpha=1$, the function $w$ belongs to $ W^{1,p}(C_{1,\theta_0})$ for any $1\le
p<{4/3}$. \fin
\noindent{\bf Proof of Lemma 3.4 : }  
\medskip
\noindent{\bf Step 1 } The norms are equivalent.

\noindent Let be $V$ in $H^1_{\rho,\Gamma_0}(\GS,\R^3)$. We applied the classical Korn
inequality   to the membrane displacements $V_M^{(l)}=V^{(l)}_{1}\vec e^{(l)}_1+ V^{(l)}_{2}\vec e^{(l)}_2$  and then  we
add all the inequalities to obtain 
$$\sum_{l=1}^N\bigl\{||\nabla V^{(l)}_1||^2_{[L^2(\omega_{l})]^2}+||\nabla V^{(l)}_2||^2_{[L^2(\omega_{l} )]^2}\bigr\}\le
C\bigl(|V|_{\rho}^2+||V||^2_{L^2(\GS,\R^3)}\bigr)$$ hence
$$||V||_{\rho}\le C\bigl\{|V|_{\rho}+||V||_{L^2(\GS,\R^3)}\bigr\}.$$ The space $H^1_{\rho,\Gamma_0}(\GS,\R^3)$ is
embedded in
$L^2(\GS,\R^3)$  (see Lemma 4.4).  Then we prove by contradiction that there exists a constant $C_0$ such
that $||V||_{\rho}\le C_0 |V|_{\rho}$. Moreover we can immediately see that there exists $C_1$ such that $|V|_{\rho}\le
C_1||V||_{\rho}.$ The norms $|\cdot|_{\rho}$ and $||\cdot||_{\rho}$ are therefore equivalent.
\vskip 1mm
\noindent{\bf Step 2  } The space $H^1_{\Gamma_0}(\GS,\R^3)\cap L^\infty(\GS,\R^3)$ is dense in
$H^1_{\rho,\Gamma_0}(\GS,\R^3)\cap L^\infty(\GS,\R^3)$.
\vskip 1mm
\noindent Let be $V\in H^1_{\rho,\Gamma_0}(\GS,\R^3)\cap L^\infty(\GS,\R^3)$, we consider the sequence
$\bigl(V_\delta\bigr)_{\delta>0}$ defined by
$$V_\delta(\widehat{x})=V(\widehat{x})\sum_{A\in {\cal N}}m\Bigl({dist(\widehat{x}, A)\over \delta}\Bigr),\qquad
\widehat{x}\in \GS$$ where $m$ is given by $(4.12)$. The displacement $V_\delta$ is equal to zero in the neighborhood of each
vertex belonging to ${\cal N}$.  This displacement  belongs to   $H^1_{\Gamma_0}(\GS,\R^3)\cap L^\infty(\GS,\R^3)$ and we
have
$$||V_\delta-V||_{L^2(\GS,\R^3)}\le C\delta||V||_{L^\infty(\GS,\R^3)},\qquad ||V_\delta||_{L^\infty(\GS,\R^3)}\le
C||V||_{L^\infty(\GS,\R^3)}.$$ We calculate the  gradient of the restriction of $V^{(l)}_\delta$  to each face
$\overline{\omega}_l$. Using  the $L^2$ estimate of $V_\delta-V$, we obtain $|V_\delta|_\rho\le
C\bigl\{|V|_\rho+||V||_{L^\infty(\GS,\R^3)}\bigr\}$.  The constant does not depend on  $\delta$. The sequence
$\bigl(V_\delta\bigr)_{\delta>0}$ weakly converges to $V$ in $H^1_{\rho,\Gamma_0}(\GS,\R^3)$, which gives the density of
$H^1_{\Gamma_0}(\GS,\R^3)\cap L^\infty(\GS,\R^3)$ into $H^1_{\rho,\Gamma_0}(\GS,\R^3)\cap L^\infty(\GS,\R^3)$.
\vskip 1mm
\noindent{\bf Step 3  } The space $H^1_{\rho,\Gamma_0}(\GS,\R^3)\cap L^\infty(\GS,\R^3)$ is dense in
$H^1_{\rho,\Gamma_0}(\GS,\R^3)$.
\vskip 1mm
\noindent We consider the truncature function $T_M$ from $\R^3$ into $\R^3$ defined by
$$ T_M(x)=\left\{\eqalign{ &x\hskip 4em\hbox{if}\; ||x||_2<M\cr &{x\over ||x||_2}M\quad\hbox{if}\; ||x||_2\ge M,\cr}\right.$$
where $M$ belongs to $\R^*_+$ and where $||\cdot||_2$ is the euclidian norm of $\R^3$. The map $T_M$ is piecewise
${\cal C}^1$  verifying $T_M(0)=0$ and $||\nabla T_M||_{[L^\infty(\R^3, \R^3)]^3}\le C$ (the constant
does not depend on $M$).

\noindent Let $V\in H^1_{\rho,\Gamma_0}(\GS,\R^3)$,  $V_M=T_M(V)$ belongs to $H^1_{\rho,\Gamma_0}(\GS,\R^3)\cap
L^\infty(\GS,\R^3)$ and verifies
$$||V_M||_{L^2(\GS,\R^3)}\le ||V||_{L^2(\GS,\R^3)},\qquad ||V_M||_{L^\infty(\GS,\R^3)}\le M,\qquad |V_M|_\rho\le
C|V|_\rho.$$ The constant does not depend on $M$. When $M$ tends to   infinity, $V_M$ tends strongly to $V$ in
$H^1_{\rho,\Gamma_0}(\GS,\R^3)$. Hence the  density of
$H^1_{ \rho,\Gamma_0}(\GS,\R^3)\cap L^\infty(\GS,\R^3)$  in $H^1_{\rho,\Gamma_0}(\GS,\R^3)$.\fin
\noindent{\bf Proof of Lemma 3.7 : } We put $V\in D_E(\GS)$. As in the proof of Lemma 3.4  we get
$$\sum_{l=1}^N\bigl\{||\nabla V^{(l)}_1||^2_{[L^2(\omega_{l})]^2}+||\nabla V^{(l)}_2||^2_{[L^2(\omega_{l} )]^2}\bigr\}\le
C\bigl(||V||_{E}^2+||V||^2_{L^2(\GS,\R^3)}\bigr)$$   We put $J\in{\cal J}$ a common edge to the faces
$\overline{\omega}_l$ and $\overline{\omega}_k$. The restrictions to $J$ of the membrane  displacements 
$V^{(l)}_M$ and $V^{(k)}_M$ completely define the restriction $V_{|J}$. Hence we get
$$\sum_{J\in{\cal J}}||V_{|J}||^2_{H^{1/ 2}(J ,\R^3)}\le C\bigl(||V||^2_{E}+||V||^2_{L^2(\GS,\R^3)}\bigr)$$
\noindent With the help of Lemma 4.4 we build a displacement $W\in H^1_{\rho,\Gamma_0}(\GS,\R^3)$
 such that
$$W_\alpha^{(l)}=V^{(l)}_\alpha,\enskip \forall l\in\{1,\ldots,N\},\qquad W_{|J}=V_{|J},\enskip \forall J\in{\cal J},$$
and verifying
$$ |W|^2_{\rho} \le C\sum_{l=1}^N\bigl\{||\nabla V^{(l)}_1||^2_{[L^2(\omega_{l})]^2}+||\nabla
V^{(l)}_2||^2_{ [L^2(\omega_l )]^2}\bigr\}+C\sum_{J\in{\cal J}}||V_{|J}||^2_{H^{1/ 2}(J ,\R^3)}\le
C\bigl(||V||^2_{E}+||V||^2_{L^2(\GS,\R^3)}\bigr).$$  The displacement  $V-W$ is of inextensional type
and hence orthogonal to $V$, hence
$$|V|_{\rho}\le |W|_{\rho}\le C\bigl(||V||_{E}+||V||_{L^2(\GS,\R^3)}\bigr).\leqno (4.16)$$ 
\noindent Now we show that the norm  $||\cdot||_{E}$ is equivalent to the norm  $|\cdot|_{\rho}$ in 
$D_{E}(\GS)$. We already have $||V||_{E}\le |V|_{\rho}$ for any $V\in D_E(\GS)$.  We suppose that the norms are not
equivalent. For any $n\in\N^*$, we can find $V_n\in D_{E}(\GS)$ such that $||V_n||_{E}\le {1/ n}$ and $|V_n|_{\rho}=1$. The
sequence   $\bigl(V_n\bigr)_{n\in\N^*}$ being bounded in $H^1_{\rho}(\GS,\R^3)$,   we can then extract  a sub-sequence, still
denoted in the same way, such that
$$V_n\rightharpoonup V\quad \hbox{weakly in}\quad H^1_{\rho,\Gamma_0}(\GS,\R^3).$$ The limit $V$
belongs also to $D_E(\GS)$. Let us make $n$ tend to infinity  in the inequality $||V_n||_E\le{1/n}$,  we obtain
$\gamma_{\alpha\beta}(V^{(l)})=0$. The displacement $V$ is of inextensional type, and hence is equal to zero. 

\noindent If $p\in[1,4/3[$ the space $H^1_{\rho}(\GS,\R^3)$ is continuously imbedded in $W^{1,p}(\GS,\R^3)$   (see Lemma
4.4) and  if $p\in ]1,4/3[$ the space $W^{1,p}(\GS,\R^3)$ is compactly imbedded in $L^2(\GS,\R^3)$.  Hence the sequence
$\bigl(V_n\bigr)_{n\in\N^*}$ converges strongly to  0 in  $L^2(\GS,\R^3)$, hence
$||V_n||_{L^2(\GS,\R^3)}\longrightarrow 0$. From $ (4.16)$ follows then that the sequence $\bigl(V_n\bigr)_{n\in\N^*}$
converges strongly to  0 in $H^1_{\rho,\Gamma_0}(\GS,\R^3)$. This stands in contradiction with $|V_n|_{\rho}=1$.\fin 
\vskip 1mm
\noindent{\bf 4.4 Annex D. The inextensional displacements}
\vskip 1mm
\noindent{\bf 4.4.a. The inextensional displacements of $D_I(\GS)$}
\vskip 1mm
\noindent Let $U$ be an inextensional displacement. From the definition of the inextensional displacements we have 
$\gamma_{\alpha\beta}(U^{(l)})=0$ in $\overline{\omega}_l$. Hence, in each face the membrane displacement
$U^{(l)}_{M}=U^{(l)}_{1}\vec e^{(l)}_1+U^{(l)}_2\vec e^{(l)}_2$ is a rigid displacement. The restriction of $U$ to an edge $J\in
{\cal J}$ is then 
$$U_{|J}(M)=\overrightarrow{A}_J+\overrightarrow{B}_J\land \overrightarrow{A_JM} \qquad \forall M\in J\leqno(4.17)$$ where
$A_J$ is an vertex of the edge. The vectors $\overrightarrow{A}_J$ and $\overrightarrow{B}_J$ depend only on the edge. We
choose  $\overrightarrow{B}_J$ orthogonal to $\vec e_J$ to have the unicity of this vector.
\vskip 1mm
\noindent{\bf 4.4.b. The inextensional displacements of ${\cal D}_I(\GS)$} 
\vskip 1mm
\noindent A displacement ${\cal A}\in {\cal D}_I(\GS)$ verifies  $\displaystyle {\partial{\cal A}^{(l)}\over
\partial x^{(l)}_i}=\widehat{\nabla}{\cal A}^{(l)}\land\vec e^{(l)}_i$. This displacement  belongs also to  $D_I(\GS)$, hence 
$\widehat{\nabla}{\cal A}^{(l)}_{3}$ is  constant in each face. Let $A$ be a vertex belonging to ${\cal N}$ and let $J$ and  $L$
be two edges sharing the vertex $A$. The functions ${\cal A}$ and $\widehat{\nabla}{\cal A}$ belong to $H^1(\GS,\R^3)$.
The hypothesis {\bf H2} implies that
$$||{\cal A}_{|J\cup L}||_{H^{1/2}(J\cup L,\R^3)}\le C||{\cal A}||_I,\qquad ||\widehat{\nabla}{\cal A}_{|J\cup
L}||_{H^{1/2}(J\cup L,\R^3)}\le C||{\cal A}||_I.\leqno(4.18)$$ From $(4.17)$ we have
${\cal A}_{|J}(M)=\overrightarrow{A}_J+\overrightarrow{B}_J\land \overrightarrow{AM}$ and
${\cal A}_{|L}(M)=\overrightarrow{A}_L+\overrightarrow{B}_L\land \overrightarrow{AM}$, hence
$\overrightarrow{A}_J=\overrightarrow{A}_L$. We denote ${\cal A}(A)$ this value which is common to all the edges containing
the vertex $A$. 

\noindent We also have $\overrightarrow{B}_J\land\vec e_J=\widehat{\nabla}{\cal A}_{|J}\land\vec e_J$. The vector
$\widehat{\nabla}{\cal A}_{|J}\land\vec e_J$ is  constant along the edge $J$, hence
$$\overrightarrow{B}_J\cdot(\vec e_J\land\vec e_L)=\overrightarrow{B}_L\cdot( \vec e_J\land\vec e_L).$$ There 
exists a vector 
$\overrightarrow{B}_{JL}\in\R^3$ (depending on $\overrightarrow{B}_J$ and $\overrightarrow{B}_L$) such that 
$$\overrightarrow{B}_J\land\vec e_J=\overrightarrow{B}_{JL}\land\vec e_J,\qquad 
\qquad\overrightarrow{B}_L\land\vec e_L=\overrightarrow{B}_{JL}\land\vec e_L.$$
\noindent Since $\widehat{\nabla}{\cal A}_{|L}(M)\land\vec e_L=\overrightarrow{B}_L\land\vec
e_L=\overrightarrow{B}_{JL}\land\vec e_L$ for any $M\in L$. From $(4.18)$ we deduce that
$$\int_0^{L_J}{|\widehat{\nabla}{\cal A}_{|J}(A+t\vec e_J)\land\vec e_L-\overrightarrow{B}_{JL}
\land\vec e_L|^2\over t}dt\le C||{\cal A}||_I,$$ hence 
$$\int_0^{L_J}{|\widehat{\nabla}{\cal A}_{|J}(A+t\vec e_J)-\overrightarrow{B}_{JL})|^2\over t}dt\le C||{\cal A}||_I,$$ because
we have  $\widehat{\nabla}{\cal A}_{|J}(M)\land\vec e_J=\overrightarrow{B}_{JL}\land\vec e_J$ for any $M\in J$. The vector
$\overrightarrow{B}_{JL}$ does not depend on the edge $L$. Hence this vector is independant from the edges that go via
$A$, and is denoted $\widehat{\nabla}{\cal A}(A)$. 

\noindent The restriction of the displacement ${\cal A}$  to any edge that goes via $A$ is the restriction to that edge of a rigid
displacement depending only on the vertex $A$,
$$\forall J\in {\cal J},\qquad\forall M\in J,\qquad {\cal A}_{|J}(M)={\cal A}(A)+\widehat{\nabla}{\cal
A}(A)\land\overrightarrow{AM}.$$
\noindent{\bf 4.5 Annex E. The test functions}
\vskip 1mm
\noindent{\bf Proof of Lemma 3.9 : } Let $U$ be an element belonging to  $H^1_{\Gamma_0}(\GS_{\delta_0},\R^3)$. 
\vskip 1mm
\noindent We suppose  that the real $\delta_0$ is such that the two balls centered in the extremities of the edges and of radius
$8\eta_0\delta_0$ do not share any common point.
\vskip 1mm
\noindent{\bf Step 1 }   For any  $\delta$ in the interval $]0,\delta_0]$, we build a displacement
$U_{\delta,1}$ constant in the neighborhood of each vertex belonging to ${\cal N}$ and approaching
$U$. 
\vskip 1mm
\noindent We begin with modifying $U$ in the neighborhood of an vertex. Let $A$ be an vertex common to the faces
$\overline{\omega}_{l_1},\; \ldots,\overline{\omega}_{l_p}$;
$\overrightarrow{A}^{(l_i)}$  the mean value of $U^{(l_i)}$ in the disc $B(A; \delta)\cap\omega_{l_i,\delta_0}$ and 
$\overrightarrow{A}$ the mean value of the vectors $\overrightarrow{A}^{(l_i)}$. If the vertex $A$ belongs to
$\Gamma_0$ we replace $\overrightarrow{A}$ by $\vec 0$.

\noindent We define the displacement $U_{\delta, A}$ in  $\overline{\omega}_{l,\delta_0}$ by
$$U^{(l)}_{\delta, A}(\widehat{ x}^{(l)})=U(\widehat{ x}^{(l)})m\Bigl({r^{(l)}\over 2\eta_0\delta}\Bigr)+\Bigl\{
1-m\Bigl({r^{(l)}\over 2\eta_0\delta}\Bigr)\Bigr\}\overrightarrow{A},\qquad l\in\{l_1,\ldots, l_p\}.$$ where $m$ has been
introduced by $(4.12)$ and where $r^{(l)}=dist(\widehat{x}^{(l)},A)$. In  $B(A;
2\eta_0\delta)\cap\overline{\omega}_{l,\delta_0}$ the displacement $U^{(l)}_{\delta,A}$ is by construction constant and equal
to $\overrightarrow{A}$.   We have 
$${\partial U^{(l)}_{\delta,A}\over \partial x_\alpha^{(l)}}-{\partial U^{(l)}\over \partial x_\alpha^{(l)}}=\Bigl\{U^{(l)}
-\overrightarrow{A}\Bigr\}m{'}\Bigl({r^{(l)}\over 2\eta_0\delta}\Bigr){x_\alpha^{(l)} \over 2\eta_0\delta
r^{(l)}}-\Bigl\{1-m\Bigl({r^{(l)}\over 2\eta_0\delta}\Bigr)\Bigr\} {\partial U^{(l)}\over \partial x^{(l)}_\alpha}.$$ Let us 
estimate the $L^2$ norm of the gradient of $U_{\delta,A}^{(l)}-U^{(l)}$, $l\in\{l_1,\ldots,l_p\}$,
$$\eqalign{ ||\nabla U^{(l)}_{\delta,A}-\nabla U^{(l)}||^2_{[L^2(\omega_{l,\delta_0},\R^3)]^2} &\le C||\nabla U^{(l)}||^2
_{[L^2(B(A;4\eta_0\delta)\cap\overline{\omega}_{l,\delta_0},\R^3)]^2} +{C\over
\delta^2}||U^{(l)}-\overrightarrow{A}||^2_{L^2(B(A;4\eta_0\delta)\cap\overline{\omega}_{l,\delta_0},\R^3)}\cr}$$
\noindent  The Poincar\'e-Wirtinger inequality allows us to estimate the   $L^2$ norm of
$U^{(l)}-\overrightarrow{A}^{(l)}$ in the disc
$B(A;\delta)\cap\omega_{l,\delta_0}$,
$$||U^{(l)}-\overrightarrow{A}^{(l)}||^2_{L^2(B(A; \delta)\cap\overline{\omega}_{l,\delta_0},\R^3)}\le C\delta^2||\nabla
U^{(l)}||^2 _{[L^2(B(A;
\delta)\cap\overline{\omega}_{l,\delta_0};\R^3)]^2}.$$ If $J$ is an edge of vertex
$A$ contained in $\overline{\omega}_l$, then we have
$$||U_{|J}-\overrightarrow{A}^{(l)}||^2_{L^2(B(A; \delta)\cap J,\R^3)}\le C\delta||\nabla U^{(l)}||^2 _{[L^2(B(A; 
\delta)\cap\overline{\omega}_{l,\delta_0},\R^3)]^2}$$ For any other face 
$\overline{\omega}_{k,\delta_0}$ containing $J\cap B(A; \delta)$, we also have the above estimate, hence
$$||\overrightarrow{A}^{(k)}-\overrightarrow{A}^{(l)}||^2_2\le C\bigl\{||\nabla U^{(l)}||^2_{[L^2(B(A; \delta)
\cap\overline{\omega}_{l,\delta_0},\R^3)]^2}+||\nabla U^{(k)}||^2_{[L^2(B(A; \delta)\cap\overline{\omega}_{k,\delta_0}
,\R^3)]^2}\bigr\}.$$
\noindent We deduce that
$$\eqalign{ ||U^{(l)}-\overrightarrow{A}||^2_{L^2(B(A; \delta)\cap\omega_{l,\delta_0},\R^3)} &\le
C\delta^2\sum_{i=1}^p||\nabla U^{(l_i)}||^2 _{[L^2(B(A;
\delta)\cap\omega_{l_i,\delta_0},\R^3)]^2}\cr
\Longrightarrow \qquad ||U^{(l)}-\overrightarrow{A}||^2_{L^2(B(A; 4\eta_0\delta)\cap\omega_{l,\delta_0},\R^3)} & \le
C\delta^2\sum_{k=1}^p||\nabla U^{(l_k)}||^2 _{[L^2(B(A; 4\eta_0\delta)\cap\omega_{l_k,\delta_0},\R^3)]^2}\cr}$$ And
eventually
$$||U_{\delta,A}-U||_{H^1(\GS_{\delta_0},\R^3)}\le C\sum_{k=1}^p||\nabla U^{(l_k)}||_{[L^2(B(A;4\eta_0\delta)\cap
\omega_{l_k,\delta_0},\R^3)]^2}.$$ We can do the same with all the structure vertexes. We denote 
$U_{\delta,1}$ the displacement obtained after having modified $U$ in a neighborhood of each vertex. Hence we have
$$||U_{\delta,1}-U||_{H^1(\GS_{\delta_0},\R^3)}\le C\sum_{A\in{\cal N}}|| U||_{H^1(B(A;4\eta_0\delta)\cap
\GS_{\delta_0},\R^3)}.$$
\noindent{\bf Step 2} Let $J$ be an edge belonging to ${\cal J}$. This edge is  common to the faces 
$\overline{\omega}_{l_1},\;  \ldots,\overline{\omega}_{l_p}$. We denote $V_{\delta,J}^{(l_k)}$ the element of
$H^1(J,\R^3)$ defined by
$$V_{\delta,J}^{(l_k)}(\widehat{x}_J)={1\over 2 \delta}\int_{- \delta}^{ \delta}U_{\delta,1}^{(l_k)}(
\widehat{x}_J+ se^{(l_k)}_{J,\perp})ds$$ and we denote $V_{\delta,J}$ the mean values of $V_{\delta,J}^{(l_k)}$
($\{e^{(l_k)}_{J},e^{(l_k)}_{J,\perp}\}$ is an orthonormal basis of the direction of the face $\overline{\omega}_{l_k}$). We have
$V_{\delta,J}^{(l_k)}\in H^1(J,\R^3)$ and
$$\left\{\eqalign{ ||U-V_{\delta,J}||_{L^2(\widehat{J}_{\delta},\R^3)}&\le C \delta||
U_{\delta,1}||_{H^1(\widehat{J}_{\delta},\R^3)}\le C \delta|| U||_{H^1(\widehat{J}_{\delta},\R^3)}\cr
||U-V_{\delta,J}||_{H^1(\widehat{J}_{\delta},\R^3)}&\le C ||U_{\delta,1}||_{H^1(\widehat{J}_{\delta},\R^3)}\le C
||U||_{H^1(\widehat{J}_{\delta},\R^3)}\cr}\right.\leqno(4.19)$$ where $\widehat{J}_{\delta}=\bigl\{\widehat{x}\in
\GS_{\delta_0}\; |\; dist(\widehat{x},J)<
\eta_0\delta\bigr\}$, $\widehat{J}_{\delta}$ is the union of two-dimensional sets of breadth $2\eta_0\delta$ and of length
 $L_J+2\eta_0\delta$. If the edge $J$ is contained in
$\Gamma_0$, we take
$V_{\delta,J}=0$, in this case we have again the estimate  $(4.19)$. The displacement
$$U_\delta=\sum_{J\in{\cal J}}U_{\delta,1}\Bigl\{1-m\Bigl({dist(\cdot, J)\over \eta_0\delta}\Bigr)\Bigr\} +\sum_{J\in{\cal
J}}V_{\delta,J}m\Bigl({dist(\cdot, J)\over \eta_0\delta}\Bigr)$$ belongs to $H^1_{\Gamma_0}({\cal S}_\delta,\R^3)$ and
verifies
$$||U_{\delta}-U||_{H^1(\GS_{\delta_0},\R^3)}\le C\sum_{A\in{\cal N}}|| U||_{H^1(B(A;4\eta_0\delta)\cap
\GS_{\delta_0},\R^3)}+ C\sum_{J\in{\cal J}}|| U||_{H^1(\widehat{J}_{\delta},\R^3)}$$
\noindent The constant is independant of $\delta$. Lemma 3.9 is proved.\fin

\noindent{\bf Proof of Lemma 3.11 : }Let be $U\in {\cal D}_I(\GS)$. We recall that there exists only one function in
$H^1_{\Gamma_0}(\GS,\R^3)$ denoted $\widehat{\nabla} U$ such that $\displaystyle {\partial U^{(l)} \over\partial
x_\alpha^{(l)}}=\widehat{\nabla}U^{(l)}\land \vec e^{(l)}_\alpha$,  $ l\in\{1,\ldots,N\}$, $\alpha\in\{1,2\}$.

\noindent{\bf Step 1 } Extension of $U^{(l)}$ and of $\widehat{\nabla} U^{(l)}$ to $\omega_{l,\delta}$. The displacement
$$u^{(l)}(x)={1\over \delta}\Bigl\{U^{(l)}(\widehat{x}^{(l)})+\widehat{\nabla}U^{(l)}(\widehat{x}^{(l)})\land x^{(l)}_3\vec 
e^{(l)}_3\Bigr\},\qquad x\in
\Omega_{l,\delta}=\omega_l\times]-\delta,\delta[,$$ of the plate $\Omega_{l,\delta}$  extends into a displacement
still denoted $u^{(l)}$ of the plate
$\Omega^{''}_{l,\delta}=\omega_{l,\delta}\times ]-\delta,\delta[$. The extension $u^{(l)}$ is by construction equal to zero on
$\Gamma_{0,\delta}\cap
\Omega^{''}_{l,\delta}$. We have
$${\cal E}(u^{(l)}, \Omega^{''}_{l,\delta})\le C{\cal E}(u^{(l)},\Omega_{l,\delta})\le 
C\delta||\widehat{\nabla}U^{(l)}||^2_{H^1(\omega_l,\R^3)}\qquad {\cal E}(u^{(l)},
\Omega^{''}_{l,\delta}\setminus\Omega_{l,\delta})\le  C\delta||\widehat{\nabla}U^{(l)}||^2_{H^1(\Gamma_{l,\delta},\R^3)}$$
where $\Gamma_{l,\delta}=\bigl\{\widehat{x}^{(l)}\in\omega_l\; |\; dist(\widehat{x}^{(l)},\partial\omega_l)<2\eta_0\delta)
\bigr\}$. Let $U_e^{(l)}$ be the e.p.d.  associated to $u^{(l)}$ by the formulas $(2.1)$, its components
$\displaystyle{1\over
\delta}{\cal U}^{(l)}$ and $\displaystyle{1\over\delta}{\cal R}^{(l)}$ are the restrictions to $\omega_{l,2\delta}$
of elements
$\bigl(\hbox{denoted}\;\displaystyle{1\over \delta}{\cal U}, \;\displaystyle{1\over\delta}{\cal R}\Bigr)$ belonging to
$H^1(\GS_{\delta},\R^3)$. They verify ${\cal U}^{(l)}_{|\omega_l}=U^{(l)}$, 
${\cal R}^{(l)}_{|\omega_l}=\widehat{\nabla}U^{(l)}$ and
$$\left\{\eqalign{
&\Bigl\|{\partial {\cal U}^{(l)}_3\over\partial x_1}+{\cal R}^{(l)}_2\Bigr\|_{L^2(\omega_{l,
\delta})}+\Bigl\|{\partial {\cal U}^{(l)}_3\over \partial x_2}-{\cal R}^{(l)}_1\Bigr\|_{L^2(\omega_{l,\delta})}
\le C\delta ||\widehat{\nabla}U^{(l)}||_{H^1(\Gamma_{l,\delta},\R^3)}\cr &\bigl\|\nabla{\cal
R}^{(l)}_1\bigr\|_{L^2(\omega_{l,\delta}\setminus
\omega_l,\R^2)}+\bigl\|\nabla{\cal R}^{(l)}_2
\bigr\|_{L^2(\omega_{l,\delta}\setminus \omega_l,\R^2)}\le C ||\widehat{\nabla}U^{(l)}||_{H^1(\Gamma_{l,\delta},\R^3)}
\cr}\right.\leqno(4.20)$$
\noindent{\bf Step 2 } We denote $V^{(l)}_\delta$ the displacement
$$V^{(l)}_\delta(x)={1\over \delta}\Bigl\{{\cal U}^{(l)}(\widehat{x}^{(l)})+{\cal R}^{(l)}(\widehat{x}^{(l)})\land
x^{(l)}_3\vec e^{(l)}_3\Bigr\},\qquad x\in
\Omega{'}_{l,\delta}.$$ We modify
$V^{(l)}_\delta$ in the neighborhood of the vertexes belonging to ${\cal N}$. 

\noindent Let $A$ be an vertex belonging to ${\cal N}$, for any edge $J$ containing $A$, we have
$$U_{|J}(x)=U(A)+\widehat{\nabla}U(A)\land \overrightarrow{Ax}=r_A(x),\qquad x\in J.$$ For any face
$\overline{\omega}_l$ containing the vertex $A$, we define the displacement $U^{(l)}_{\delta, A}$ by
$$U^{(l)}_{\delta, A}(x)=V^{(l)}_\delta(x)m\Bigl({dist(x,A)\over 2\eta_0\delta}\Bigr)+{1\over \delta}\Bigl\{
1-m\Bigl({dist(x,A)\over 2\eta_0\delta}\Bigr)\Bigr\}r_A(x),\qquad x\in \Omega{'}_{l,\delta},$$ where $m$ has been
introduced by $(4.7)$. We have $U^{(l)}_{\delta, A}(x)={1/\delta}r_A(x)$ in $B(A,2\eta_0\delta)\cap
\Omega{'}_{l,\delta}$. Let us  remind that for any edge $J$ containing the vertex $A$, we have
${\cal R}_{|J}=\widehat{\nabla}U_{|J}$ and ${\cal R}_{|J}\land\vec e_J=\widehat{\nabla}U_{|J}\land\vec
e_J=\widehat{\nabla}U(A)\land\vec e_J$. Hence, thanks to  $(4.20)$, we have the following inequality :
$$\eqalign{ &||{\cal R}-\widehat{\nabla}U(A)||_{L^2(B(A,4\eta_0\delta)\cap \GS_\delta,\R^3)}+ \Bigl\|{\partial {\cal
U}^{(l)}_3\over \partial x_\alpha}-\widehat{\nabla}U(A)\land\vec e^{(l)}_\alpha\Bigr\|_{L^2(B(A,4\eta_0\delta)\cap
\GS_\delta,\R^3)}\cr
\le &C\delta\bigl\{||\widehat{\nabla}U||_{ H^1(B(A,4\eta_0\delta)\cap \GS_\delta,\R^3)}+ ||\widehat{\nabla}U^{(l)}||_{
H^1(\Gamma_{l,\delta},\R^3)}\bigr\}.\cr}$$  This inequality implies that
$$\eqalign{ ||U^{(l)}_{\delta, A}-V^{(l)}_\delta||^2_{L^2(B(A,4\eta_0\delta)\cap \Omega{'}_{l,\delta},\R^3)} &\le
C\delta^3\bigl\{||\widehat{\nabla}U||_{ H^1(B(A,4\eta_0\delta)\cap \GS_\delta,\R^3)}+ ||\widehat{\nabla}U^{(l)}||_{
H^1(\Gamma_{l,\delta},\R^3)}\bigr\},\cr
\hbox{and}\qquad{\cal E}(U^{(l)}_{\delta, A}-V^{(l)}_\delta,\Omega{'}_{l,\delta}) &\le C\delta\bigl\{ ||\widehat{\nabla}U||_{
H^1(B(A,4\eta_0\delta)\cap \GS_\delta,\R^3)}+ ||\widehat{\nabla}U^{(l)}||_{ H^1(\Gamma_{l,\delta},\R^3)}\bigr\}.\cr}$$
The displacement 
$\displaystyle U^{(l)}_{\delta, 1}=\sum_{A\in \overline{\omega}_l\cap {\cal N}}U^{(l)}_{\delta, A}$ coincides with a rigid
displacement independent of $l$ in the neighborhood of each vertex contained in the face $\overline{\omega}_l$ and verifies
$${\cal E}(U^{(l)}_{\delta, 1}-V^{(l)}_\delta,\Omega{'}_{l,\delta}) \le C\delta\bigl\{\sum_{A\in
\overline{\omega}_l\cap {\cal N}}||\widehat{\nabla}U||_{ H^1(B(A,4\eta_0\delta)\cap \GS_\delta,\R^3)}+
||\widehat{\nabla}U^{(l)}||_{ H^1(\Gamma_{l,\delta},\R^3)}\bigr\}.$$
\noindent{\bf Step 3} We modify $U^{(l)}_{\delta, 1}$ in the neighborhood of each edge belonging to ${\cal J}$.

\noindent Let $J$ be an edge belonging to several faces and $\overline{\omega}_l$ a face of $\GS$ containing
$J$. We take the orthornormal frame $(O_J;\vec e_J, \vec e^{(l)}_{J,\perp}, \vec e^{(l)}_3)$
linked to the edge $J$ and to the plate
$\Omega_{l,\delta}$ containing this edge,  $O_J$ is an extremity of $J$ and  $\vec e_{J}$  the direction of the edge
($x^{(l)}_{1J}=\widehat{x}^{(l)}\cdot \vec e_J\in [0,L_J]$, $L_J$ the length of the edge). In this frame we consider the 
neighborhood of $J$
$$J^{(l)}_{\eta\delta}=]0,L_J[\times]-\eta\delta,\eta\delta [\times ]-\delta,\delta[,\qquad \eta\ge 1.$$  The
restriction of $U^{(l)}_{\delta, 1}$ to
$J^{(l)}_\delta$ ($\eta=1$) is decomposed into the sum of an elementary rod displacement
 $U^{(l)}_{e,J}$ and of a residual displacement,
$$U^{(l)}_{e,J}(x_{1J}, x^{(l)}_{2J}, x^{(l)}_{3})={\cal U}^{(l)}_{J}(x_{1J})+{\cal
R}^{(l)}_{J}(x_{1J})\land\bigl(x^{(l)}_{2J}\vec e^{(l)}_{J,\perp} +x^{(l)}_3\vec e^{(l)}_3\bigr).$$  Let us remind that (see [4])
the components ${\cal U}^{(l)}_J$ and ${\cal R}^{(l)}_J$ of $U^{(l)}_{e,J}$  belong to 
$H^1(J,\R^3)$, and verify
$$\left\{\eqalign{ &\delta^2\Bigl\|{d{\cal R}^{(l)}_{J}\over dx_{1J}}\Bigr\|^2_{L^2(J,\R^3)}+ \Bigl\|{d{\cal U}^{(l)}_{J}\over
dx_{1J}} -{\cal R}^{(l)}_{J}\land\vec e_J\Bigr\|^2_{L^2(J,\R^3)}\le {C\over \delta^2}{\cal E}(U^{(l)}_{\delta,
1},J^{(l)}_\delta)\le {C\over \delta^2}{\cal E}(u^{(l)}, J^{(l)}_\delta)\cr  &{\cal
D}(U^{(l)}_{\delta,1}-U^{(l)}_{e,J},J^{(l)}_\delta)+{1\over\delta^2}||U^{(l)}_{\delta,1}-U^{(l)}_{e,J}||^2_{L^2(J^{(l)}
_\delta,\R^3)}\le C{\cal E}(U^{(l)}_{\delta, 1},J^{(l)}_\delta)\le C{\cal E}(u^{(l)},J^{(l)}_\delta)\cr}\right.$$ By construction,
the displacement $U^{(l)}_{e,J}$ coincides with a rigid displacement in the neighborhood of the edge extremities. We deduce that
$$||U_{|J}-{\cal U}^{(l)}_{J}||^2_{L^2(J,\R^3)}\le  C{\cal E}( u^{(l)},J^{(l)}_\delta),\qquad ||\widehat{\nabla}U_{|J}-{\cal
R}^{(l)}_{J}||^2_{L^2(J,\R^3)}\le  {C\over\delta^2}{\cal E}( u^{(l)},J^{(l)}_\delta)$$ The edge $J$ belongs to the faces
$\overline{\omega}_{l_1},\ldots , \overline{\omega}_{l_p}$. Let $U_{e,J}$ be the elementary  rod displacement equal to the
mean value of the displacements $U^{(l_1)}_{e,J},\ldots , U^{(l_p)}_{e,J}$. The components of  $U_{e,J}$ being
${\cal U}_J$ and ${\cal R}_J$, we have
$$\left\{\eqalign{ &\delta^2\Bigl\|{d{\cal R}_{J}\over dx_{1J}}\Bigr\|^2_{L^2(J,\R^3)}+ \Bigl\|{d{\cal U}_{J}\over dx_{1J}}
-{\cal R}_{J}\land\vec e_J\Bigr\|^2_{L^2(J,\R^3)}\le {C\over \delta^2}\sum_{i=1}^p{\cal E}(u^{(l_i)}, J^{(l_i)}_\delta)\cr 
&\sum_{i=1}^p\Bigl\{{\cal D}(U^{(l_i)}_{\delta,1}-U_{e,J},J^{(l_i)}_\delta)+{1\over\delta^2}||U^{(l_i)}_{\delta,1}
-U_{e,J}||^2_{ L^2(J^{(l_i)}_\delta,\R^3)}\Bigr\}\le  C\sum_{i=1}^p{\cal E}(u^{(l_i)},J^{(l_i)}_\delta)\cr}\right.$$  We
deduce  (see [4]) that
$$\left\{\eqalign{ &\sum_{i=1}^p{\cal E}(U_{e,J},J^{(l_i)}_{2\eta_0\delta})\le C \sum_{i=1}^p{\cal E}(u^{(l_i)},
J^{(l_i)}_\delta)\cr  &\sum_{i=1}^p\Bigl\{{\cal
E}(U^{(l_i)}_{\delta,1}-U_{e,J},J^{(l_i)}_{2\eta_0\delta})+{1\over\delta^2}||U^{(l_i)}_{\delta,1} -U_{e,J}||^2_{
L^2(J^{(l_i)}_{2\eta_0\delta},\R^3)}\Bigr\}\le  C\sum_{i=1}^p{\cal E}(u^{(l_i)},J^{(l_i)}_{2\eta_0\delta})\cr}\right.$$Now we
modify the displacement $U^{(l)}_{\delta,1}$ in the neighborhood of the edge $J$,
$$U^{(l)}_{\delta, J}(x)=U^{(l)}_{\delta, 1}(x)\Bigl\{1-m\Bigl({dist(\widehat{x}^{(l)},J)\over \eta_0\delta}\Bigr)\Bigr\}+
U_{e,J}(x)m\Bigl({dist(\widehat{x}^{(l)},J)\over \eta_0\delta}\Bigr)\Bigr\},\qquad x\in \Omega{'}_{l,\delta}.$$ Then we have
$$\sum_{i=1}^p{\cal E}(U^{(l_i)}_{\delta, J}-U^{(l_i)}_{\delta,1}, \Omega{'}_{l_i,\delta})\le  C\sum_{i=1}^p{\cal
E}(u^{(l_i)},J^{(l_i)}_{ 2\eta_0\delta})$$ Eventually the displacement $U_\delta$ obtained by modifying
$U^{(l)}_{\delta,1}$ in the neighborhood of each edge of ${\cal J}$ belongs to
$H^1_{\Gamma_0}({\cal S}_\delta,\R^3)$ and for any $l\in\{1,\ldots,N\}$ verifies $(3.23)$.\fin
\bigskip
\noindent{\bf REFERENCES }
\bigskip
\noindent [1] P.G. Ciarlet, Ph. Destuynder. A justification of two-dimensional linear plate model, J. M\'ecanique, Vol. 18 (2.2),
 1979, 315--344.

\noindent [2] D. Cioranescu, A. Damlamian and  G. Griso, Periodic Unfolding and Homogenization. CRAS, Ser. I 335 (2002) 99-104.

\noindent [3] Ph. Destuynder. Une th\'eorie asymptotique des plaques minces en \'elasticit\'e lin\'eaire. RMA 2, Masson (1986).

\noindent [4] G. Griso.  Asymptotic behavior of  curved rods by the unfolding method. (To appear).

\noindent [5] G. Griso.  Asymptotic behavior of structures made of curved rods. (To appear).

\noindent [6] G. Griso. Asymptotic behavior of structures made of plates.  CRAS, Ser. I 336 (2003) 101-106.

\noindent [7] H. Le Dret. Modeling of a folded plate, Comput.  Mech., 5 (1990), 401--416.

\bye